\documentclass{amsart}
\usepackage{amsmath,amscd,amssymb,amsthm,amsbsy,amscd,stmaryrd}
\usepackage{mathrsfs}
\input{xy}

\usepackage{pstricks}
\usepackage{color}
\usepackage{graphicx}
\usepackage[all]{xy}
\usepackage {hyperref}

\newtheorem{theorem}{Theorem}[section]

\newtheorem{proposition}[theorem]{Proposition}
\newtheorem{corollary}[theorem]{Corollary}
\newtheorem{remark}[theorem]{Remark}
\newtheorem{question}[theorem]{Question}
\newtheorem{conjecture}[theorem]{Conjecture}
\theoremstyle{definition}

\numberwithin{equation}{section}

\begin{document}
\title{Mabuchi and Aubin-Yau functionals over complex
three-folds}
\author{Yi Li}
\address{Department of Mathematics, Harvard University, Cambridge, MA 02138}

\email{yili@math.harvard.edu}

\begin{abstract} In this paper we construct Mabuchi
$\mathcal{L}^{{\rm M}}_{\omega}$ functional and Aubin-Yau
functionals $\mathcal{I}^{{\rm AY}}_{\omega}, \mathcal{J}^{{\rm
AY}}_{\omega}$ on any compact complex three-folds. The method
presented here will be used in the forthcoming paper \cite{L1} on
the construction of those functionals on any compact complex
manifolds, which generalizes the previous work \cite{L}.
\end{abstract}
\maketitle

\tableofcontents
\section{Introduction}

Mabuchi and Aubin-Yau functionals play a crucial role in studying
K\"ahler-Einstein metrics and constant scalar curvatures. How to
generalize these functionals from K\"ahler geometry to complex
geometry is an interesting problem. In \cite{L}, the author solved
this problem in dimension two and proved similar results in the
K\"ahler setting. By carefully checking and using a trick, we can
construct those functionals on any compact complex three-folds.
Moreover, the idea in this paper will be used in the forthcoming
paper in which we deal with higher dimension cases.

\subsection{Mabuchi and Aubin-Yau functionals on K\"ahler manifolds}
Let $(X,\omega)$ be a compact K\"ahler manifold of the complex
dimension $n$. then the volume
\begin{equation}
V_{\omega}:=\int_{X}\omega^{n}\label{1.1}
\end{equation}
depends only on the K\"ahler class of $\omega$. Let
$\mathcal{P}_{\omega}$ denote the space of K\"ahler potentials and
define the Mabuchi functional, for any smooth functions
$\varphi',\varphi''\in\mathcal{P}_{\omega}$, by
\begin{equation}
\mathcal{L}^{{\rm M},{\rm
Kahler}}_{\omega}(\varphi',\varphi''):=\frac{1}{V_{\omega}}\int^{1}_{0}\int_{X}\dot{\varphi}_{t}
\omega^{n}_{\varphi_{t}}dt\label{1.2}
\end{equation}
where $\varphi_{t}$ is any smooth path in $\mathcal{P}_{\omega}$
from $\varphi'$ to $\varphi''$. Mabuchi \cite{M} showed that
(\ref{1.2}) is well-defined.

Using (\ref{1.2}) we can define Aubin-Yau functionals as follows:
\begin{eqnarray}
\mathcal{I}^{{\rm AY},{\rm
Kahler}}_{\omega}(\varphi)&=&\frac{1}{V_{\omega}}\int_{Z}\varphi(\omega^{n}-\omega^{n}_{\varphi}),
\label{1.3}\\
\mathcal{J}^{{\rm AY},{\rm
Kahler}}_{\omega}(\varphi)&=&-\mathcal{L}^{{\rm M},{\rm
Kahler}}_{\omega}(0,\varphi)+\frac{1}{V_{\omega}}\int_{X}\varphi\omega^{n}\label{1.4}.
\end{eqnarray}
So Aubin-Yau functionals are also well-defined.

However, if $\omega$ is not closed, then the above definitions do
not make sense. Hence we should add some extra terms on the
definitions of those functionals. These extra terms should involve
$\partial\omega$ and $\overline{\partial}\omega$, but, the essential
question is to find the structure of the extra terms. In the next
section, we will answer this question.

\subsection{Mabuchi and Aubin-Yau functionals on complex
three-folds}
Throughout the rest part of this paper, we denote by $(X,g)$ a
compact complex manifold of the complex dimension $3$, and $\omega$
the associated real $(1,1)$-form. Let
\begin{equation}
\mathcal{P}_{\omega}:=\{\varphi\in C^{\infty}(X)_{\mathbb{R}}|
\omega_{\varphi}:=\omega+\sqrt{-1}\partial\overline{\partial}\varphi>0\}\label{1.5}
\end{equation}
be the space of all real-valued smooth functions on $X$ whose
associated real $(1,1)$-forms are positive.

For any $\varphi',\varphi''\in\mathcal{P}_{\omega}$, we define
\begin{eqnarray}
\mathcal{L}^{{\rm
M}}_{\omega}(\varphi',\varphi'')&=&\frac{1}{V_{\omega}}\int^{1}_{0}\int_{X}\dot{\varphi}_{t}\omega^{3}_{\varphi_{t}}dt
\nonumber\\
&-&\frac{1}{V_{\omega}}\int^{1}_{0}\int_{X}3\sqrt{-1}\partial\omega\wedge(\overline{\partial}\dot{\varphi}_{t}\cdot\varphi_{t})\wedge\omega_{\varphi_{t}}dt
\nonumber\\
&+&\frac{1}{V_{\omega}}\int^{1}_{0}\int_{X}3\sqrt{-1}\overline{\partial}\omega\wedge(\partial\dot{\varphi}_{t}\cdot\varphi_{t})\wedge\omega_{\varphi_{t}}dt
\label{1.6}\\
&-&\frac{1}{V_{\omega}}\int^{1}_{0}\int_{X}\partial\varphi_{t}\wedge\overline{\partial}\varphi_{t}\wedge\partial\omega\wedge
\overline{\partial}\dot{\varphi}_{t}-\frac{1}{V_{\omega}}\int^{1}_{0}\int_{X}\overline{\partial}\varphi_{t}\wedge\partial\varphi_{t}\wedge\overline{\partial}
\omega\wedge\partial\dot{\varphi}_{t}\nonumber
\end{eqnarray}
where $\{\varphi_{t}\}_{0\leq t\leq1}$ is any smooth path in
$\mathcal{P}_{\omega}$ from $\varphi'$ to $\varphi''$. Our first
result is

\begin{theorem} \label{thm1.1} The functional (\ref{1.6}) is independent of the choice of
the smooth path $\{\varphi_{t}\}_{0\leq t\le 1}$ in
$\mathcal{P}_{\omega}$.
\end{theorem}

For any $\varphi\in\mathcal{P}_{\omega}$ we set
\begin{equation}
\mathcal{L}^{{\rm M}}_{\omega}(\varphi):=\mathcal{L}^{{\rm
M}}_{\omega}(0,\varphi).\label{1.7}
\end{equation}
Then we have an explicit formula of $\mathcal{L}^{{\rm
M}}_{\omega}(\varphi)$:

\begin{corollary} \label{cor1.2} One has
\begin{eqnarray}
\mathcal{L}^{{\rm M}}_{\omega}(\varphi)&=&
\frac{1}{4V_{\omega}}\sum^{3}_{i=0}\int_{X}\varphi\omega^{i}_{\varphi}\wedge\omega^{3-i}
\nonumber\\
&-&\sum^{1}_{i=0}\frac{i+1}{2V_{\omega}}\int_{X}\varphi\omega^{i}_{\varphi}\wedge\omega^{1-i}\wedge\sqrt{-1}\partial\omega\wedge\overline{\partial}\varphi
\label{1.8}\\
&+&\sum^{1}_{i=0}\frac{i+1}{2V_{\omega}}\int_{X}\varphi\omega^{i}_{\varphi}\wedge\omega^{1-i}\wedge\sqrt{-1}\overline{\partial}\omega\wedge\partial\varphi.\nonumber
\end{eqnarray}
\end{corollary}

Now we define Aubin-Yau functionals $\mathcal{I}^{{\rm
AY}}_{\omega}, \mathcal{J}^{{\rm AY}}_{\omega}$ for any compact
complex three-fold $(X,\omega)$:
\begin{eqnarray}
& & \mathcal{I}^{{\rm AY}}_{\omega}(\varphi) \ \ = \ \
\frac{1}{V_{\omega}}\int_{X}\varphi(\omega^{3}-\omega^{3}_{\varphi})\label{1.9}\\
&-&\frac{3}{2V_{\omega}} \int_{X}\varphi\omega_{\varphi}
\wedge\sqrt{-1}\partial\omega\wedge\overline{\partial}\varphi-
\frac{3}{2V_{\omega}}\int_{X}\varphi\omega\wedge\sqrt{-1}\partial\omega\wedge
\overline{\partial}\varphi \nonumber\\
&+&\frac{3}{2V_{\omega}}\int_{X}\varphi\omega_{\varphi}
\wedge\sqrt{-1}\overline{\partial}\omega\wedge\partial\varphi
+\frac{3}{2V_{\omega}}\int_{X}\varphi\omega\wedge\sqrt{-1}\overline{\partial}\omega\wedge\partial\varphi
\nonumber,\\
& & \mathcal{J}^{{\rm AY}}_{\omega}(\varphi) \ \ = \ \
-\mathcal{L}^{{\rm
M}}_{\omega}(\varphi)+\frac{1}{V_{\omega}}\int_{X}\varphi\omega^{3}
\label{1.10}\\
&-&\frac{3}{2V_{\omega}}\int_{X}\varphi\omega_{\varphi}
\wedge\sqrt{-1}\partial\omega\wedge\overline{\partial}\varphi
-\frac{3}{2V_{\omega}}\int_{X}\varphi\omega\wedge\sqrt{-1}\partial\omega\wedge\overline{\partial}\varphi\nonumber\\
&+&\frac{3}{2V_{\omega}}\int_{X}\varphi\omega_{\varphi}
\wedge\sqrt{-1}\overline{\partial}\omega\wedge\partial\varphi
+\frac{3}{2V_{\omega}}\int_{X}\varphi\omega\wedge\sqrt{-1}\overline{\partial}\omega\wedge\partial\varphi.\nonumber
\end{eqnarray}

An important result is

\begin{theorem} \label{1.3} For any $\varphi\in\mathcal{P}_{\omega}$, one has
\begin{eqnarray}
\frac{3}{4}\mathcal{I}^{{\rm
AY}}_{\omega}(\varphi)-\mathcal{J}^{{\rm
AY}}_{\omega}(\varphi)&\geq&0, \label{1.11}\\
4\mathcal{J}^{{\rm AY}}_{\omega}(\varphi)-\mathcal{I}^{{\rm
AY}}_{\omega}(\varphi)&\geq&0.\label{1.12}
\end{eqnarray}
In particular
\begin{eqnarray*}
\frac{1}{4}\mathcal{I}^{{\rm
AY}}_{\omega}(\varphi)&\leq&\mathcal{J}^{{\rm AY}}_{\omega}(\varphi)
\ \ \leq \ \ \frac{3}{4}\mathcal{I}^{{\rm AY}}_{\omega}(\varphi),
\\
\frac{4}{3}\mathcal{J}^{{\rm
AY}}_{\omega}(\varphi)&\leq&\mathcal{I}^{{\rm AY}}_{\omega}(\varphi)
\ \ \leq \ \ 4\mathcal{J}^{{\rm AY}}_{\omega}(\varphi), \\
\frac{1}{3}\mathcal{J}^{{\rm AY}}_{\omega}(\varphi) \ \ \leq \ \
\frac{1}{4}\mathcal{J}^{{\rm AY}}_{\omega}(\varphi)&\leq&
\mathcal{I}^{{\rm
AY}}_{\omega}(\varphi)-\mathcal{J}^{{\rm AY}}_{\omega}(\varphi)\\
&\leq&\frac{3}{4}\mathcal{I}^{{\rm AY}}_{\omega}(\varphi) \ \ \leq \
\ 3\mathcal{J}^{{\rm AY}}_{\omega}(\varphi).
\end{eqnarray*}
\end{theorem}

\subsection{Volume estimate}
For any compact K\"ahler manifold the volume (\ref{1.1}) depends
only on the K\"ahler class, but, this fact doesn't hold for the
general compact Hermitian manifolds. To understand the change of
volumes for compact Hermitian three-fold $(X,\omega)$, the author
 \cite{L} introduced three quantities associated to $\omega$:
\begin{eqnarray}
{\rm
Err}_{\omega}(\varphi)&:=&\int_{X}\omega^{n}-\int_{X}\omega^{n}_{\varphi},
\ \ \ \varphi\in\mathcal{P}_{\omega},
\label{1.13}\\
{\rm
SupErr}_{\omega}&:=&\sup_{\varphi\in\mathcal{P}^{0}_{\omega}}({\rm
Err}_{\omega}(\varphi)), \label{1.14}\\
{\rm
InfErr}_{\omega}&:=&\inf_{\varphi\in\mathcal{P}^{0}_{\omega}}({\rm
Err}_{\omega}(\varphi)).\label{1.15}
\end{eqnarray}
where
$\mathcal{P}^{0}_{\omega}=\{\varphi\in\mathcal{P}_{\omega}|\sup_{X}\varphi=0\}$
is the normalized subspace of $\mathcal{P}_{\omega}$. In any case,
we have
\begin{equation}
{\rm InfErr}_{\omega}\leq0\leq{\rm
SupErr}_{\omega}\leq\int_{X}\omega^{n}.\label{1.16}
\end{equation}
It implies that the quantity ${\rm SupErr}_{\omega}$ is bounded, but
the quantity ${\rm InfErr}_{\omega}$ may not have a lower bound.
Notice that the lower boundedness of ${\rm InfErr}_{\omega}$ is
equivalent to the upper boundedness of
$\sup_{\varphi\in\mathcal{P}^{0}_{\omega}}\int_{X}\omega^{n}_{\varphi}$.
We can ask a natural question \cite{LV}:

\begin{question} \label{q1.4} Under what condition (weaker than the K\"ahler
condition), the quantity
\begin{equation}
\sup_{\varphi\in\mathcal{P}^{0}_{\omega}}\int_{X}\omega^{n}_{\varphi}
\label{1.17}
\end{equation}
is bounded from above?
\end{question}

\begin{remark} \label{rmk1.5} (1) In \cite{L} the author showed that if
$\partial\overline{\partial}\omega=\partial\omega\wedge\overline{\partial}\omega=0$,
then the volume $\int_{X}\omega^{n}_{\varphi}$ does not depend on
the choice of the smooth function $\varphi\in\mathcal{P}_{\omega}$.
Consequently, $\int_{X}\omega^{n}_{\varphi}=\int_{X}\omega^{n}$ for
any function $\varphi\in\mathcal{P}_{\omega}$, which gives an
affirmative answer to Question \ref{q1.4}. In the surface case, we
can only assume $\partial\overline{\partial}\omega=0$. \\
(3) Let $(X,\omega)$ be a compact Hermitian manifold of the complex
dimension $2$. By a theorem of Ganduchon \cite{G}, there exists a
smooth function $u$, unique up to a constant, such that
\begin{equation}
\partial\overline{\partial}\omega_{{\rm G}}=0, \ \ \
\omega_{{\rm G}}:=e^{u}\omega.
\end{equation}
V. Tosatti and B. Weinkove \cite{TW1} showed that (\ref{1.17}) is
bounded by $\int_{X}(2e^{u-\inf_{X}u}-1)\omega^{2}$ from above.
\end{remark}

Using the idea in \cite{TW1} we can show that, under a suitable
condition, ${\rm InfErr}_{\omega}$ is bounded from below when $n=3$
and, consequently, (\ref{1.17}) has an upper bound.

\begin{theorem} \label{thm1.6} Suppose that $(X,g)$ is a compact Hermitian manifold
of the complex dimension $3$ and $\omega$ is its associated real
$(1,1)$-form. If $\partial\overline{\partial}\omega=0$, then ${\rm
InfErr}_{\omega}$ is bounded from below. More precisely, we have
\begin{equation}
{\rm InfErr}_{\omega}\geq3\left(1-e^{2\cdot{\rm
osc}(u)}\right)\cdot\int_{X}\omega^{3}.\label{1.19}
\end{equation}
Here $u$ is a real-valued smooth function on $X$ such that
$\omega_{G}=e^{u}\cdot \omega$ is a Gauduchon metric, i.e.,
$\partial\overline{\partial}(\omega^{2}_{G})=0$. In particular
\begin{equation}
\int_{X}\omega^{3}\leq\sup_{\varphi\in\mathcal{P}^{0}_{\omega}}\int_{X}\omega^{3}_{\varphi}\leq
\left(3e^{2\cdot{\rm
osc}(u)}-2\right)\int_{X}\omega^{3}.\label{1.20}
\end{equation}
\end{theorem}

Another interesting question \cite{LV} is

\begin{conjecture} \label{c1.7} For any compact Hermitian manifold $(X,\omega)$,
one has
\begin{equation}
\inf_{\varphi\in\mathcal{P}^{0}_{\omega}}\int_{X}\omega^{n}_{\varphi}>0.\label{1.21}
\end{equation}
\end{conjecture}

\begin{remark} \label{rmk1.8} (1) When $\omega$ is K\"ahler, for any smooth
function $\varphi\in\mathcal{P}_{\omega}$, we have
$\int_{X}\omega^{n}_{\varphi}=\int_{X}\omega^{n}$. Hence
$\inf_{\varphi\in\mathcal{P}^{0}_{\omega}}\int_{X}\omega^{n}_{\varphi}=\int_{X}\omega^{n}>0$.\\
(2) In \cite{TW1}, the authors confirmed Conjecture \ref{c1.7} for
$n=2$ provided that $u$ satisfies $\sup_{X}(u)-\inf_{X}(u)<{\rm
log}2$. More precisely
\begin{equation}
\inf_{\varphi\in\mathcal{P}^{0}_{\omega}}\int_{X}\omega^{n}_{\varphi}\geq\int_{X}\left(
2e^{u-\sup_{X}(u)}-1\right)\omega^{2}>0.\label{1.22}
\end{equation}
\end{remark}

We can prove Conjecture \ref{c1.7} for $n=3$ provided that a similar
condition holds for $u$.

\begin{theorem} \label{thm1.9} Suppose that $(X,g)$ is a compact Hermitian manifold
of the complex dimension $3$ and $\omega$ is its associated real
$(1,1)$-form. We select a real-valued smooth function $u$ on $X$ so
that $e^{u}\cdot \omega$ is a Gauduchon metric. If
\begin{equation*}
{\rm osc}(u)=\sup_{X}(u)-\inf_{X}(u)\leq\frac{1}{2}\cdot{\rm
ln}\frac{3}{2}, \ \ \
\partial\overline{\partial}\omega=0,
\end{equation*}
then
\begin{equation}
\inf_{\varphi\in\mathcal{P}^{0}_{\omega}}\int_{X}\omega^{3}_{\varphi}\geq\int_{X}\omega^{3}>0.\label{1.23}
\end{equation}
\end{theorem}

\section{Mabuchi $\mathcal{L}^{{\rm M}}_{\omega}$ functional on
compact complex three-folds}
Let $(X,g)$ be a compact complex manifold of the complex dimension
$3$, and $\omega$ be its associated real $(1,1)$-form. In this
section we define Mabuchi $\mathcal{L}^{{\rm M}}_{\omega}$
functional and prove the independence of the choice of the smooth
path. As a consequence, we give an explicit formula for
$\mathcal{L}^{{\rm M}}_{\omega}$.

Let $\varphi',\varphi''\in\mathcal{P}_{\omega}$ and
$\{\varphi_{t}\}_{0\leq t\leq1}$ be a smooth path in
$\mathcal{P}_{\omega}$ from $\varphi'$ to $\varphi''$. We
define
\begin{equation}
\mathcal{L}^{0}_{\omega}(\varphi',\varphi'')
=\frac{1}{V_{\omega}}\int^{1}_{0}\int_{X}\dot{\varphi}_{t}
\omega^{3}_{\varphi_{t}}dt.\label{2.1}
\end{equation}
Set
\begin{equation}
\psi(s,t):=s\cdot\varphi_{t}, \ \ \ 0\leq s,t\leq1.\label{2.2}
\end{equation}
Consider a $1$-form on $[0,1]\times[0,1]$
\begin{equation}
\Psi^{0}:=\left(\int_{X}\frac{\partial\psi}{\partial
s}\cdot\omega^{3}_{\psi}\right)ds+\left(\int_{X}\frac{\partial\psi}{\partial
t}\cdot\omega^{3}_{\psi}\right)dt.\label{2.3}
\end{equation}
Taking the differential on $\Psi^{0}$, we have
\begin{equation*}
d\Psi^{0}=I^{0}\cdot dt\wedge ds
\end{equation*}
where the quantity $I^{0}$ is given by
\begin{equation}
I^{0}=\int_{X}\frac{\partial}{\partial
t}\left(\frac{\partial\psi}{\partial
s}\cdot\omega^{3}_{\psi}\right)-\int_{X}\frac{\partial}{\partial
s}\left(\frac{\partial\psi}{\partial
t}\cdot\omega^{3}_{\psi}\right).\label{2.4}
\end{equation}
As in \cite{L} we simplify $I^{0}$ in two slightly different ways.
Directly computing shows
\begin{eqnarray*}
I^{0}&=&\int_{X}\left[\frac{\partial^{2}\psi}{\partial t\partial
s}\cdot\omega^{3}_{\psi}+\frac{\partial\psi}{\partial
s}\cdot3\omega^{2}_{\psi}\wedge\sqrt{-1}\partial\overline{\partial}\left(\frac{\partial
\psi}{\partial t}\right)\right] \\
&-&\int_{X}\left[\frac{\partial^{2}\psi}{\partial s\partial
t}\cdot\omega^{3}_{\psi}+\frac{\partial\psi}{\partial
t}\wedge3\omega^{2}_{\psi}\wedge\sqrt{-1}\partial\overline{\partial}\left(\frac{\partial\psi}{\partial
s}\right)\right] \\
&=&\int_{X}3\frac{\partial\psi}{\partial
s}\omega^{2}_{\psi}\wedge\sqrt{-1}\partial\overline{\partial}\left(\frac{\partial\psi}{\partial
t}\right)+\int_{X}3\frac{\partial\psi}{\partial
t}\omega^{2}_{\psi}\wedge\sqrt{-1}\overline{\partial}\partial\left(\frac{\partial\psi}{\partial
s}\right) \\
&=&-\int_{X}3\frac{\partial\psi}{\partial
s}\omega^{2}_{\psi}\wedge\sqrt{-1}\overline{\partial}\partial\left(\frac{\partial\psi}{\partial
t}\right)-\int_{X}3\frac{\partial\psi}{\partial
t}\omega^{2}_{\psi}\wedge\sqrt{-1}\partial\overline{\partial}\left(\frac{\partial\psi}{\partial
s}\right).
\end{eqnarray*}
Notice that the last two steps are the differential expressions for
$I^{0}$. Using the first expression, we have
\begin{eqnarray}
I^{0}&=&\int_{X}3\frac{\partial\psi}{\partial
s}\omega^{2}_{\psi}\wedge\sqrt{-1}\partial\overline{\partial}\left(\frac{\partial\psi}{\partial
t}\right)+\int_{X}3\frac{\partial\psi}{\partial
t}\omega^{2}_{\psi}\wedge\sqrt{-1}\overline{\partial}\partial\left(\frac{\partial\psi}{\partial
s}\right)\nonumber\\
&=&\int_{X}-3\sqrt{-1}\partial\left(\frac{\partial\psi}{\partial
s}\omega^{2}_{\psi}\right)\wedge\overline{\partial}\left(\frac{\partial\psi}{\partial
t}\right)+\int_{X}-3\sqrt{-1}\overline{\partial}\left(\frac{\partial\psi}{\partial
t}\omega^{2}_{\psi}\right)\wedge\partial\left(\frac{\partial\psi}{\partial
s}\right)\nonumber \\
&=&\int_{X}-3\sqrt{-1}\left[\partial\left(\frac{\partial\psi}{\partial
s}\right)\wedge\omega^{2}_{\psi}+\frac{\partial\psi}{\partial
s}2\omega_{\psi}\wedge\partial\omega\right]\wedge\overline{\partial}\left(\frac{\partial\psi}{\partial
t}\right) \nonumber\\
&=&\int_{X}-3\sqrt{-1}\left[\overline{\partial}\left(\frac{\partial\psi}{\partial
t}\right)\wedge\omega^{2}_{\psi}+\frac{\partial\psi}{\partial
t}2\omega_{\psi}\wedge\overline{\partial}\omega\right]\wedge\partial\left(\frac{\partial\psi}{\partial
s}\right) \nonumber\\
&=&\int_{X}-6\sqrt{-1}\frac{\partial\psi}{\partial
s}\omega_{\psi}\wedge\partial\omega\wedge\overline{\partial}\left(\frac{\partial\psi}{\partial
t}\right)+\int_{X}-6\sqrt{-1}\frac{\partial\psi}{\partial
t}\omega_{\psi}\wedge\overline{\partial}\omega\wedge\partial\left(\frac{\partial\psi}{\partial
s}\right)\label{2.5}.
\end{eqnarray}
Similarly, we have
\begin{eqnarray}
I^{0}&=&-\int_{X}3\frac{\partial\psi}{\partial
s}\omega^{2}_{\psi}\wedge\sqrt{-1}\overline{\partial}\partial\left(\frac{\partial\psi}{\partial
t}\right)-\int_{X}3\frac{\partial\psi}{\partial
t}\omega^{2}_{\psi}\wedge\sqrt{-1}\partial\overline{\partial}\left(\frac{\partial\psi}{\partial
s}\right) \nonumber\\
&=&\int_{X}6\sqrt{-1}\frac{\partial\psi}{\partial
s}\omega_{\psi}\wedge\overline{\partial}\omega\wedge\partial\left(\frac{\partial\psi}{\partial
t}\right)+\int_{X}6\sqrt{-1}\frac{\partial\psi}{\partial
t}\omega_{\psi}\wedge\partial\omega\wedge\overline{\partial}\left(\frac{\partial\psi}{\partial
s}\right)\label{2.6}.
\end{eqnarray}
Therefore, from (\ref{2.5}) and (\ref{2.6}) it follows that
\begin{eqnarray}
\frac{2I^{0}}{6\sqrt{-1}}&=&\int_{X}\overline{\partial}\left(\frac{\partial\psi}{\partial
t}\right)\frac{\partial\psi}{\partial
s}\wedge\omega_{\psi}\wedge\partial\omega+\int_{X}\partial\left(\frac{\partial\psi}{\partial
s}\right)\frac{\partial\psi}{\partial
t}\wedge\omega_{\psi}\wedge\overline{\partial}\omega \nonumber \\
&-&\int_{X}\overline{\partial}\left(\frac{\partial\psi}{\partial
s}\right)\frac{\partial\psi}{\partial
t}\wedge\omega_{\psi}\wedge\partial\omega-\int_{X}\partial\left(\frac{\partial\psi}{\partial
t}\right)\frac{\partial\psi}{\partial
s}\wedge\omega_{\psi}\wedge\overline{\partial}\omega.\label{2.7}
\end{eqnarray}

Next we define
\begin{eqnarray}
\mathcal{L}^{1}_{\omega}(\varphi',\varphi'')&=&
\frac{1}{V_{\omega}}\int^{1}_{0}\int_{X}a_{1}\cdot\partial\omega\wedge
\omega_{\varphi_{t}}\wedge(\overline{\partial}\dot{\varphi}_{t}\cdot\varphi_{t})dt,
\label{2.8}\\
\mathcal{L}^{2}_{\omega}(\varphi',\varphi'')&=&\frac{1}{V_{\omega}}\int^{1}_{0}\int_{X}
a_{2}\cdot\overline{\partial}\omega\wedge\omega_{\varphi_{t}}
\wedge(\partial\dot{\varphi}_{t}\cdot\varphi_{t})dt.\label{2.9}
\end{eqnarray}
Here $a_{1},a_{2}$ are non-zero constants and we require
$\overline{a_{1}}=a_{2}$. Similarly, we consider
\begin{eqnarray*}
\Psi^{1}&=&\left[\int_{X}a_{1}\partial\omega\wedge\omega_{\psi}\wedge\left(
\overline{\partial}\left(\frac{\partial\psi}{\partial
s}\right)\psi\right)\right]ds+\left[\int_{X}a_{1}\partial\omega\wedge\omega_{\psi}\wedge\left(
\overline{\partial}\left(\frac{\partial\psi}{\partial
t}\right)\psi\right)\right]dt, \\
\Psi^{2}&=&\left[\int_{X}a_{2}\overline{\partial}\omega\wedge\omega_{\psi}\wedge\left(\partial\left(\frac{\partial\psi}{\partial
s}\right)\psi\right)\right]ds+\left[\int_{X}a_{2}\overline{\partial}\omega\wedge\omega_{\psi}\wedge
\left(\partial\left(\frac{\partial\psi}{\partial
t}\right)\psi\right)\right]dt.
\end{eqnarray*}
Therefore
\begin{equation*}
d\Psi^{1}=I^{1}\cdot dt\wedge ds
\end{equation*}
where
\begin{equation*}
I^{1}=\int_{X}a_{1}\frac{\partial}{\partial
t}\left[\partial\omega\wedge\omega_{\psi}\wedge\left(\overline{\partial}\left(\frac{\partial\psi}{\partial
s}\right)\cdot\psi\right)\right]-
\int_{X}a_{1}\frac{\partial}{\partial
s}\left[\partial\omega\wedge\omega_{\psi}\wedge\left(\overline{\partial}\left(\frac{\partial\psi}{\partial
t}\right)\cdot\psi\right)\right].
\end{equation*}
Consequently,
\begin{eqnarray*}
\frac{I^{1}}{a_{1}}&=&\int_{X}-\frac{\partial}{\partial
t}\left[\left(\psi\cdot\overline{\partial}\left(\frac{\partial\psi}{\partial
s}\right)\right)\wedge\omega_{\psi}\wedge\partial\omega\right]+
\int_{X}\frac{\partial}{\partial
s}\left[\left(\psi\cdot\overline{\partial}\left(\frac{\partial\psi}{\partial
t}\right)\right)\wedge\omega_{\psi}\wedge\partial\omega\right] \\
&=&\int_{X}-\left[\frac{\partial\psi}{\partial
t}\cdot\overline{\partial}\left(\frac{\partial\psi}{\partial
s}\right)+\psi\cdot\overline{\partial}\left(\frac{\partial^{2}\psi}{\partial
t\partial s}\right)\right]\wedge\omega_{\psi}\wedge\partial\omega \\
&+&\int_{X}\left[\frac{\partial\psi}{\partial
s}\cdot\overline{\partial}\left(\frac{\partial\psi}{\partial
t}\right)+\psi\cdot\overline{\partial}\left(\frac{\partial^{2}\psi}{\partial
s\partial t}\right)\right]\wedge\omega_{\psi}\wedge\partial\omega \\
&+&\int_{X}\psi\cdot\overline{\partial}\left(\frac{\partial\psi}{\partial
s}\right)\wedge-\sqrt{-1}\partial\overline{\partial}\left(\frac{\partial\psi}{\partial
t}\right)\wedge\partial\omega \\
&+&\int_{X}\psi\cdot\overline{\partial}\left(\frac{\partial\psi}{\partial
t}\right)\wedge\sqrt{-1}\partial\overline{\partial}\left(\frac{\partial\psi}{\partial
s}\right)\wedge\partial\omega \\
&=&\int_{X}-\frac{\partial\psi}{\partial
t}\cdot\overline{\partial}\left(\frac{\partial\psi}{\partial
s}\right)\wedge\omega_{\psi}\wedge\partial\omega+\int_{X}\frac{\partial\psi}{\partial
s}\cdot\overline{\partial}\left(\frac{\partial\psi}{\partial
t}\right)\wedge\omega_{\psi}\wedge\partial\omega \\
&+&\int_{X}\psi\cdot\overline{\partial}\left(\frac{\partial\psi}{\partial
s}\right)\wedge-\sqrt{-1}\partial\overline{\partial}\left(\frac{\partial\psi}{\partial
t}\right)\wedge\partial\omega \\
&+&\int_{X}\psi\cdot\overline{\partial}\left(\frac{\partial\psi}{\partial
t}\right)\wedge\sqrt{-1}\partial\overline{\partial}\left(\frac{\partial\psi}{\partial
s}\right)\wedge\partial\omega.
\end{eqnarray*}

In the same fashion way, we have
\begin{equation*}
d\Psi^{2}=I^{2}\cdot dt\wedge ds,
\end{equation*}
where
\begin{eqnarray*}
\frac{I^{2}}{a_{2}}&=&\int_{X}-\frac{\partial\psi}{\partial
t}\cdot\partial\left(\frac{\partial\psi}{\partial
s}\right)\wedge\omega_{\psi}\wedge\overline{\partial}\omega
+\int_{X}\frac{\partial\psi}{\partial
s}\cdot\partial\left(\frac{\partial\psi}{\partial
t}\right)\wedge\omega_{\psi}\wedge\overline{\partial}\omega \\
&+&\int_{X}\psi\cdot\partial\left(\frac{\partial\psi}{\partial
s}\right)\wedge\sqrt{-1}\overline{\partial}\partial\left(\frac{\partial\psi}{\partial
t}\right)\wedge\overline{\partial}\omega \\
&+&\int_{X}\psi\cdot\partial\left(\frac{\partial\psi}{\partial
t}\right)\wedge-\sqrt{-1}\overline{\partial}\partial\left(\frac{\partial\psi}{\partial
s}\right)\wedge\overline{\partial}\omega.
\end{eqnarray*}
To simplify notation we set
\begin{eqnarray}
\mathcal{A}&:=&\int_{X}\psi\overline{\partial}\left(\frac{\partial\psi}{\partial
s}\right)\wedge\partial\omega\wedge-\sqrt{-1}\partial\overline{\partial}\left(\frac{\partial\psi}{\partial
t}\right),\label{2.10}\\
\mathcal{B}&:=&\int_{X}\psi\overline{\partial}\left(\frac{\partial\psi}{\partial
t}\right)\wedge\partial\omega\wedge\sqrt{-1}\partial\overline{\partial}\left(\frac{\partial\psi}{\partial
s}\right).\label{2.11}
\end{eqnarray}
Using the above symbols gives
\begin{eqnarray*}
\frac{I^{1}}{a_{1}}&=&\int_{X}-\frac{\partial\psi}{\partial
t}\cdot\overline{\partial}\left(\frac{\partial\psi}{\partial
s}\right)\wedge\omega_{\psi}\wedge\partial\omega+\int_{X}\frac{\partial\psi}{\partial
s}\cdot\overline{\partial}\left(\frac{\partial\psi}{\partial
t}\right)\wedge\omega_{\psi}\wedge\partial\omega
+\mathcal{A}+\mathcal{B}, \\
\frac{I^{2}}{a_{2}}&=&\int_{X}-\frac{\partial\psi}{\partial
t}\cdot\partial\left(\frac{\partial\psi}{\partial
s}\right)\wedge\omega_{\psi}\wedge\overline{\partial}\omega
+\int_{X}\frac{\partial\psi}{\partial
s}\cdot\partial\left(\frac{\partial\psi}{\partial
t}\right)\wedge\omega_{\psi}\wedge\overline{\partial}\omega+
\overline{\mathcal{A}}+\overline{\mathcal{B}},
\end{eqnarray*}
and the last two terms can be determined completely as follows:
\begin{eqnarray*}
\mathcal{A}&=&\int_{X}\sqrt{-1}\partial\left[\psi\overline{\partial}\left(\frac{\partial\psi}{\partial
s}\right)\wedge\partial\omega\right]\wedge\overline{\partial}\left(\frac{\partial\psi}{\partial
t}\right) \\
&=&\int_{X}\sqrt{-1}\left[\partial\left(\psi\cdot\overline{\partial}\left(\frac{\partial\psi}{\partial
s}\right)\right)\wedge\partial\omega\right]\wedge\overline{\partial}\left(\frac{\partial\psi}{\partial
t}\right) \\
&=&\int_{X}\sqrt{-1}\left[\partial\psi\wedge\overline{\partial}\left(\frac{\partial\psi}{\partial
s}\right)+\psi\cdot\partial\overline{\partial}\left(\frac{\partial\psi}{\partial
s}\right)\right]\wedge\partial\omega\wedge\overline{\partial}\left(\frac{\partial\psi}{\partial
t}\right)\\
&=&\int_{X}-\psi\overline{\partial}\left(\frac{\partial\psi}{\partial
t}\right)\wedge\partial\omega\wedge\sqrt{-1}\partial\overline{\partial}\left(\frac{\partial\psi}{\partial
s}\right)\\
&+&\int_{X}-\sqrt{-1}\partial\psi\wedge\partial\omega\wedge\overline{\partial}\left(\frac{\partial\psi}{\partial
s}\right)\wedge\overline{\partial}\left(\frac{\partial\psi}{\partial
t}\right) \\
&=&-\mathcal{B}+\int_{X}-\sqrt{-1}\partial\psi\wedge\partial\omega\wedge\overline{\partial}\left(\frac{\partial\psi}{\partial
s}\right)\wedge\overline{\partial}\left(\frac{\partial\psi}{\partial
t}\right).
\end{eqnarray*}
Adding the term $\mathcal{B}$ on both sides we obtain
\begin{eqnarray}
\mathcal{A}+\mathcal{B}&=&\int_{X}-\sqrt{-1}\partial\psi\wedge\partial\omega\wedge\overline{\partial}
\left(\frac{\partial\psi}{\partial
s}\right)\wedge\overline{\partial}\left(\frac{\partial\psi}{\partial
t}\right), \label{2.12}\\
\overline{\mathcal{A}}+\overline{\mathcal{B}}&=&\int_{X}\sqrt{-1}\overline{\partial}\psi
\wedge\overline{\partial}\omega\wedge\partial\left(\frac{\partial\psi}{\partial
s}\right)\wedge\partial\left(\frac{\partial\psi}{\partial
t}\right).\label{2.13}
\end{eqnarray}

The final step is to introduce
\begin{eqnarray}
\mathcal{L}^{3}_{\omega}(\varphi',\varphi'')&:=&\frac{1}{V_{\omega}}
\int^{1}_{0}\int_{X}a_{3}\partial\varphi_{t}\wedge\partial\omega\wedge\overline{\partial}\dot{\varphi}_{t}
\wedge\overline{\partial}\varphi_{t}, \label{2.14}\\
\mathcal{L}^{4}_{\omega}(\varphi',\varphi'')&:=&\frac{1}{V_{\omega}}
\int^{1}_{0}\int_{X}a_{4}\overline{\partial}\varphi_{t}\wedge\overline{\partial}\omega\wedge\partial
\dot{\varphi}_{t}\wedge\partial\varphi_{t},\label{2.15}
\end{eqnarray}
where $a_{3},a_{4}$ are non-zero constants determined later and we
require $\overline{a_{3}}=a_{4}$. Consider
\begin{eqnarray*}
\Psi^{3}&=&\left[\int_{X}a_{3}\partial\psi\wedge\partial\omega\wedge\overline{\partial}\left(\frac{\partial\psi}{\partial
s}\right)\wedge\overline{\partial}\psi\right]ds\\
&+&\left[\int_{X}a_{3}\partial\psi\wedge\partial\omega\wedge\overline{\partial}\left(\frac{\partial\psi}{\partial
t}\right)\wedge\overline{\partial}\psi\right]dt, \\
\Psi^{4}&=&\left[\int_{X}a_{4}\overline{\partial}\psi\wedge\overline{\partial}\omega\wedge
\partial\left(\frac{\partial\psi}{\partial
s}\right)\wedge\partial\psi\right]ds\\
&+&\left[\int_{X}a_{4}\overline{\partial}\psi\wedge\overline{\partial}\omega\wedge\partial
\left(\frac{\partial\psi}{\partial
t}\right)\wedge\partial\psi\right]dt.
\end{eqnarray*}
We take the differential on $\Psi^{3}$ and $\Psi^{4}$, and these
differentials can be written as
\begin{equation*}
d\Psi^{3}=I^{3}\cdot dt\wedge ds, \ \ \ d\Psi^{4}=I^{4}\cdot
dt\wedge ds
\end{equation*}
where
\begin{eqnarray*}
I^{3}&=&\int_{X}a_{3}\frac{\partial}{\partial
t}\left[\partial\psi\wedge\partial\omega\wedge\overline{\partial}\left(\frac{\partial\psi}{\partial
s}\right)\wedge\overline{\partial}\psi\right]\\
&-&\int_{X}a_{3}\frac{\partial}{\partial
s}\left[\partial\psi\wedge\partial\omega\wedge\overline{\partial}\left(\frac{\partial\psi}{\partial
t}\right)\wedge\overline{\partial}\psi\right], \\
I^{4}&=&\int_{X}a_{4}\frac{\partial}{\partial
t}\left[\overline{\partial}\psi\wedge\overline{\partial}\omega\wedge\partial\left(\frac{\partial\psi}{\partial
s}\right)\wedge\partial\psi\right]\\
&-&\int_{X}a_{4}\frac{\partial}{\partial
s}\left[\overline{\partial}\psi\wedge\overline{\partial}\omega\wedge\partial\left(\frac{\partial\psi}{\partial
t}\right)\wedge\partial\psi\right].
\end{eqnarray*}
Calculate
\begin{eqnarray*}
\frac{I^{3}}{a_{3}}&=&\int_{X}\left[\partial\left(\frac{\partial\psi}{\partial
t}\right)\wedge\partial\omega\wedge\overline{\partial}\left(\frac{\partial\psi}{\partial
s}\right)\wedge\overline{\partial}\psi+\partial\psi\wedge\partial\omega\wedge\overline{\partial}
\left(\frac{\partial^{2}\psi}{\partial t\partial
s}\right)\wedge\overline{\partial}\psi\right.\\
&+&\left.\partial\psi\wedge\partial\omega\wedge\overline{\partial}\left(\frac{\partial\psi}{\partial
s}\right)\wedge\overline{\partial}\left(\frac{\partial\psi}{\partial
t}\right)\right] \\
&-&\int_{X}\left[\partial\left(\frac{\partial\psi}{\partial
s}\right)\wedge\partial\omega\wedge\overline{\partial}\left(\frac{\partial\psi}{\partial
t}\right)\wedge\overline{\partial}\psi+\partial\psi\wedge\partial\omega\wedge\overline{\partial}
\left(\frac{\partial^{2}\psi}{\partial s\partial
t}\right)\wedge\overline{\partial}\psi\right. \\
&+&\left.\partial\psi\wedge\partial\omega\wedge\overline{\partial}\left(\frac{\partial\psi}{\partial
t}\right)\wedge
\overline{\partial}\left(\frac{\partial\psi}{\partial
s}\right)\right] \\
&=&\int_{X}-\partial\left(\frac{\partial\psi}{\partial
t}\right)\wedge\overline{\partial}\left(\frac{\partial\psi}{\partial
s}\right)\wedge\partial\omega\wedge\overline{\partial}\psi+\int_{X}\partial\left(\frac{\partial\psi}{\partial
s}\right)\wedge\overline{\partial}\left(\frac{\partial\psi}{\partial
t}\right)\wedge\partial\omega\wedge\overline{\partial}\psi \\
&+&2\int_{X}\partial\psi\wedge\partial\omega\wedge\overline{\partial}\left(\frac{\partial\psi}{\partial
s}\right)\wedge\overline{\partial}\left(\frac{\partial\psi}{\partial
t}\right) \\
&=&\int_{X}-\partial\left(\frac{\partial\psi}{\partial
t}\right)\wedge\overline{\partial}\left(\frac{\partial\psi}{\partial
s}\right)\wedge\partial\omega\wedge\overline{\partial}\psi+\int_{X}\partial\left(\frac{\partial\psi}{\partial
s}\right)\wedge\overline{\partial}\left(\frac{\partial\psi}{\partial
t}\right)\wedge\partial\omega\wedge\overline{\partial}\psi \\
&+&\frac{2}{-\sqrt{-1}}(\mathcal{A}+\mathcal{B}).
\end{eqnarray*}
Similarly,
\begin{eqnarray*}
\frac{I^{4}}{a_{4}}&=&\int_{X}-\overline{\partial}\left(\frac{\partial\psi}{\partial
t}\right)\wedge\partial\left(\frac{\partial\psi}{\partial
s}\right)\wedge\overline{\partial}\omega\wedge\partial\psi+\int_{X}\overline{\partial}
\left(\frac{\partial\psi}{\partial
s}\right)\wedge\partial\left(\frac{\partial\psi}{\partial
t}\right)\wedge\overline{\partial}\omega\wedge\partial\psi \\
&+&\frac{2}{\sqrt{-1}}(\overline{\mathcal{A}}+\overline{\mathcal{B}}).
\end{eqnarray*}
If we set
\begin{equation*}
\mathcal{H}:=\int_{X}-\partial\left(\frac{\partial\psi}{\partial
t}\right)\wedge\overline{\partial}\left(\frac{\partial\psi}{\partial
s}\right)\wedge\partial\omega\wedge\overline{\partial}\psi+\int_{X}\partial\left(\frac{\partial\psi}{\partial
s}\right)\wedge\overline{\partial}\left(\frac{\partial\psi}{\partial
t}\right)\wedge\partial\omega\wedge\overline{\partial}\psi,
\end{equation*}
then those two terms have the shorted expressions:
\begin{equation*}
\frac{I^{3}}{a_{3}}=\mathcal{H}+\frac{2}{-\sqrt{-1}}(\mathcal{A}+\mathcal{B}),
\ \ \
\frac{I^{4}}{a_{4}}=\overline{\mathcal{H}}+\frac{2}{\sqrt{-1}}(\overline{\mathcal{A}}+\overline{\mathcal{B}}).
\end{equation*}
On the other hand, directly by definition, we have
\begin{eqnarray*}
\mathcal{H}&=&\int_{X}\left[\overline{\partial}\left(\frac{\partial\psi}{\partial
s}\right)\wedge\partial\omega\wedge\overline{\partial}\psi\right]\wedge\partial\left(\frac{\partial\psi}{\partial
t}\right)+\int_{X}\left[\partial\left(\frac{\partial\psi}{\partial
s}\right)\wedge\partial\omega\wedge\overline{\partial}\psi\right]
\wedge\overline{\partial}\left(\frac{\partial\psi}{\partial
t}\right) \\
&=&\int_{X}\partial\left[\overline{\partial}\left(\frac{\partial\psi}{\partial
s}\right)\wedge\partial\omega\wedge\overline{\partial}\psi\right]\frac{\partial\psi}{\partial
t}+\int_{X}\overline{\partial}\left[\partial\left(\frac{\partial\psi}{\partial
s}\right)\wedge\partial\omega\wedge\overline{\partial}\psi\right]\frac{\partial\psi}{\partial
t} \\
&=&\int_{X}\frac{\partial\psi}{\partial
t}\left[\partial\overline{\partial}\left(\frac{\partial\psi}{\partial
s}\right)\wedge\partial\omega\wedge\overline{\partial}\psi+\overline{\partial}\left(\frac{\partial\psi}{\partial
s}\right)\wedge\partial\omega\wedge\partial\overline{\partial}\psi\right]
\\
&+&\int_{X}\frac{\partial\psi}{\partial
t}\left[\overline{\partial}\partial\left(\frac{\partial\psi}{\partial
s}\right)\wedge\partial\omega\wedge\overline{\partial}\psi-\partial\left(\frac{\partial\psi}{\partial
s}\right)\wedge\overline{\partial}\partial\omega\wedge\overline{\partial}\psi\right]
\\
&=&\int_{X}\frac{\partial\psi}{\partial
t}\cdot\overline{\partial}\left(\frac{\partial\psi}{\partial
s}\right)\wedge\partial\omega\wedge\partial\overline{\partial}\psi-\int_{X}\frac{\partial\psi}{\partial
t}\cdot\partial\left(\frac{\partial\psi}{\partial
s}\right)\wedge\overline{\partial}\partial\omega\wedge\overline{\partial}\psi
\\
&=&\int_{X}\overline{\partial}\left[\frac{\partial\psi}{\partial
t}\cdot\overline{\partial}\left(\frac{\partial\psi}{\partial
s}\right)\wedge\partial\omega\right]\wedge\partial\psi-\int_{X}\frac{\partial\psi}{\partial
t}\cdot\partial\left(\frac{\partial\psi}{\partial
s}\right)\wedge\overline{\partial}\partial\omega\wedge\overline{\partial}\psi
\\
&=&\int_{X}\overline{\partial}\left(\frac{\partial\psi}{\partial
t}\right)\wedge\overline{\partial}\left(\frac{\partial\psi}{\partial
s}\right)\wedge\partial\omega\wedge\partial\psi-\int_{X}\frac{\partial\psi}{\partial
t}\cdot\overline{\partial}\left(\frac{\partial\psi}{\partial
s}\right)\wedge\overline{\partial}\partial\omega\wedge\partial\psi
\\
&-&\int_{X}\frac{\partial\psi}{\partial
t}\cdot\partial\left(\frac{\partial\psi}{\partial
s}\right)\wedge\overline{\partial}\partial\omega\wedge\overline{\partial}\psi
\\
&=&\int_{X}\partial\psi\wedge\partial\omega\wedge\overline{\partial}\left(\frac{\partial\psi}{\partial
s}\right)\wedge\overline{\partial}\left(\frac{\partial\psi}{\partial
t}\right) \\
&+&\int_{X}\frac{\partial\psi}{\partial
t}\cdot\overline{\partial}\left(\frac{\partial\psi}{\partial
s}\right)\wedge\partial\overline{\partial}\omega\wedge\partial\psi+\int_{X}
\frac{\partial\psi}{\partial
t}\cdot\partial\left(\frac{\partial\psi}{\partial
s}\right)\wedge\partial\overline{\partial}\omega\wedge\overline{\partial}\psi
\end{eqnarray*}
Taking the complex conjugate gives
\begin{eqnarray*}
\overline{\mathcal{H}}&=&\int_{X}\overline{\partial}\psi\wedge\overline{\partial}\omega
\wedge\partial\left(\frac{\partial\psi}{\partial
s}\right)\wedge\partial\left(\frac{\partial\psi}{\partial t}\right)
\\
&-&\int_{X}\frac{\partial\psi}{\partial
t}\cdot\partial\left(\frac{\partial\psi}{\partial
s}\right)\wedge\partial\overline{\partial}\omega\wedge\overline{\partial}\psi-\int_{X}
\frac{\partial\psi}{\partial
t}\cdot\overline{\partial}\left(\frac{\partial\psi}{\partial
s}\right)\wedge\partial\overline{\partial}\omega\wedge\partial\psi.
\end{eqnarray*}
Hence
\begin{eqnarray*}
\mathcal{H}+\overline{\mathcal{H}}&=&\int_{X}\partial\psi\wedge\partial\omega\wedge\overline{\partial}\left(\frac{\partial\psi}{\partial
s}\right)\wedge\overline{\partial}\left(\frac{\partial\psi}{\partial
t}\right)+\int_{X}\overline{\partial}\psi\wedge\overline{\partial}\omega
\wedge\partial\left(\frac{\partial\psi}{\partial
s}\right)\wedge\partial\left(\frac{\partial\psi}{\partial t}\right)
\\
&=&\frac{\mathcal{A}+\mathcal{B}}{-\sqrt{-1}}+\frac{\overline{\mathcal{A}}+\overline{\mathcal{B}}}{\sqrt{-1}}
\ \ = \ \
\frac{(\mathcal{A}+\mathcal{B})-(\overline{\mathcal{A}}+\overline{\mathcal{B}})}{-\sqrt{-1}}.
\end{eqnarray*}
Consequently
\begin{eqnarray*}
\frac{2I^{0}}{6\sqrt{-1}}&=&\frac{I^{1}}{a_{1}}-\frac{I^{2}}{a_{2}}+(\mathcal{A}
+\mathcal{B})-(\overline{\mathcal{A}}+\overline{\mathcal{B}}), \\
\frac{I^{3}}{a_{3}}+\frac{I^{4}}{a_{4}}&=&\mathcal{H}+\overline{\mathcal{H}}
+\frac{2}{-\sqrt{-1}}\left[(\mathcal{A}+\mathcal{B})-(\overline{\mathcal{A}}+\overline{\mathcal{B}})\right]
\ \ = \ \
3\sqrt{-1}[(\mathcal{A}+\mathcal{B})-(\overline{\mathcal{A}}+\overline{\mathcal{B}})]
\end{eqnarray*}
it follows that
\begin{equation*}
\frac{I^{0}}{3\sqrt{-1}}=\frac{I^{0}}{a_{1}}-\frac{I^{2}}{a_{2}}+
\frac{1}{3\sqrt{-1}}\left(\frac{I^{3}}{a_{3}}+\frac{I^{4}}{a_{4}}\right).
\end{equation*}
or
\begin{equation*}
I^{0}=\frac{3\sqrt{-1}}{a_{1}}I^{1}-\frac{3\sqrt{-1}}{a_{2}}I^{2}
+\frac{I^{3}}{a_{3}}+\frac{I^{4}}{a_{4}}.
\end{equation*}
Selecting
\begin{equation*}
a_{1}=-3\sqrt{-1}, \ \ \ a_{2}=3\sqrt{-1}, \ \ \ a_{3}=a_{4}=-1
\end{equation*}
we deduce
\begin{equation}
I^{0}+I^{1}+I^{2}+I^{3}+I^{4}=0.\label{2.16}
\end{equation}

\begin{theorem} \label{thm2.1} The functional
\begin{eqnarray}
\mathcal{L}^{{\rm
M}}_{\omega}(\varphi',\varphi'')&=&\frac{1}{V_{\omega}}\int^{1}_{0}\int_{X}\dot{\varphi}_{t}\omega^{3}_{\varphi_{t}}dt
\nonumber\\
&-&\frac{1}{V_{\omega}}\int^{1}_{0}\int_{X}3\sqrt{-1}\partial\omega\wedge(\overline{\partial}\dot{\varphi}_{t}\cdot\varphi_{t})\wedge\omega_{\varphi_{t}}dt
\nonumber\\
&+&\frac{1}{V_{\omega}}\int^{1}_{0}\int_{X}3\sqrt{-1}\overline{\partial}\omega\wedge(\partial\dot{\varphi}_{t}\cdot\varphi_{t})\wedge\omega_{\varphi_{t}}dt
\label{2.17}\\
&-&\frac{1}{V_{\omega}}\int^{1}_{0}\int_{X}\partial\varphi_{t}\wedge\overline{\partial}\varphi_{t}\wedge\partial\omega\wedge
\overline{\partial}\dot{\varphi}_{t}-\frac{1}{V_{\omega}}\int^{1}_{0}\int_{X}\overline{\partial}\varphi_{t}\wedge\partial\varphi_{t}\wedge\overline{\partial}
\omega\wedge\partial\dot{\varphi}_{t}\nonumber
\end{eqnarray}
is independent of the choice of the smooth path
$\{\varphi_{t}\}_{0\leq t\leq 1}$ in $\mathcal{P}_{\omega}$, where
$\varphi_{0}=\varphi'$ and $\varphi_{1}=\varphi''$.
\end{theorem}

\begin{proof} It immediately follows from (\ref{2.16}).
\end{proof}

\begin{corollary} \label{cor2.2} For any $\varphi\in\mathcal{P}_{\omega}$ one has
\begin{eqnarray}
\mathcal{L}^{{\rm M}}_{\omega}(\varphi)&:=&\mathcal{L}^{{\rm
M}}_{\omega}(0,\varphi) \ \ = \ \
\frac{1}{4V_{\omega}}\sum^{3}_{i=0}\int_{X}\varphi\omega^{i}_{\varphi}\wedge\omega^{3-i}
\nonumber\\
&-&\sum^{1}_{i=0}\frac{i+1}{2V_{\omega}}\int_{X}\varphi\omega^{i}_{\varphi}\wedge\omega^{1-i}\wedge\sqrt{-1}\partial\omega\wedge\overline{\partial}\varphi
\label{2.11}\\
&+&\sum^{1}_{i=0}\frac{i+1}{2V_{\omega}}\int_{X}\varphi\omega^{i}_{\varphi}\wedge\omega^{1-i}\wedge\sqrt{-1}\overline{\partial}\omega\wedge\partial\varphi.\nonumber
\end{eqnarray}
\end{corollary}

\begin{proof} In Theorem \ref{2.1} we take
$\varphi_{t}=t\cdot\varphi$, then the last two terms vanish and
\begin{eqnarray*}
\mathcal{L}^{{\rm
M}}_{\omega}(\varphi)&=&\frac{1}{V_{\omega}}\int^{1}_{0}\int_{X}\varphi\omega^{3}_{t\varphi}dt
-\frac{1}{V_{\omega}}\int^{1}_{0}\int_{X}3\sqrt{-1}\partial\omega\wedge(\overline{\partial}\varphi\cdot
t\varphi)\wedge\omega_{t\varphi}dt \\
&+&\frac{1}{V_{\omega}}\int^{1}_{0}\int_{X}3\sqrt{-1}\overline{\partial}\omega\wedge(\partial\varphi\cdot
t\varphi)\wedge\omega_{t\varphi}dt:=J_{0}+J_{1}+J_{2}.
\end{eqnarray*}
Three terms are computed as follows by using the identity
$\omega_{t\varphi}=t\omega_{\varphi}+(1-t)\omega$:
\begin{eqnarray*}
J_{0}&=&\frac{1}{V_{\omega}}\int^{1}_{0}\int_{X}\varphi[t\omega_{\varphi}+(1-t)\omega]^{3}dt
\\
&=&\frac{1}{V_{\omega}}\int^{1}_{0}\int_{X}\varphi\sum^{3}_{i=0}\omega^{i}_{\varphi}\wedge\omega^{3-i}\binom{3}{i}t^{i}(1-t)^{3-i}dt
\\
&=&\frac{1}{V_{\omega}}\sum^{3}_{i=0}\int_{X}\varphi\omega^{i}_{\varphi}\wedge\omega^{3-i}\cdot
\int^{1}_{0}\binom{3}{i}t^{i}(1-t)^{3-i}dt
=\frac{1}{4V_{\omega}}\sum^{3}_{i=0}\int_{X}\varphi\omega^{i}_{\varphi}\wedge\omega^{3-i}.
\end{eqnarray*}
For $J_{1}$, we have
\begin{eqnarray*}
J_{1}&=&-\frac{1}{V_{\omega}}\int^{1}_{0}\int_{X}3\sqrt{-1}\partial\omega\wedge(\overline{\partial}\varphi\cdot
t\varphi)\wedge\omega_{t\varphi}dt \\
&=&-\frac{1}{V_{\omega}}\int^{1}_{0}\int_{X}3\sqrt{-1}\partial\omega\wedge\overline{\partial}\varphi\cdot
t\varphi\wedge[t\omega_{\varphi}+(1-t)\omega]dt \\
&=&-\frac{1}{V_{\omega}}\int_{X}3\sqrt{-1}\varphi\cdot\partial\omega\wedge\overline{\partial}\varphi\wedge\omega_{\varphi}\cdot\int^{1}_{0}t^{2}dt\\
&-&\frac{1}{V_{\omega}}\int_{X}3\sqrt{-1}\varphi\cdot\partial\omega\wedge\overline{\partial}\varphi\wedge\omega\cdot\int^{1}_{0}t(1-t)dt
\\
&=&-\frac{1}{V_{\omega}}\int_{X}\sqrt{-1}\varphi\cdot\partial\omega\wedge\overline{\partial}\varphi\wedge\omega_{\varphi}-\frac{1}{V_{\omega}}\int_{X}\frac{\sqrt{-1}}{2}\varphi\cdot\partial\omega\wedge\overline{\partial}\varphi\wedge\omega.
\end{eqnarray*}
Taking the complex conjugate gives the third term $J_{2}$.
\end{proof}

\begin{remark} \label{rmk2.3} (1) When $(X,g)$ is a compact K\"ahler three-fold, the
functional (\ref{2.1}) or (\ref{2.11}) coincides with the original
one.\\
(2) The last two terms in (\ref{2.17}) may not be zero since
$\overline{\partial}\varphi_{t}\wedge\overline{\partial}\dot{\varphi}_{t}$
is not identically zero in general. For instance, take
$\varphi_{t}=t\varphi''+(1-t)\varphi'$; then
\begin{eqnarray*}
\overline{\partial}\varphi_{t}\wedge\overline{\partial}\dot{\varphi}_{t}&=&\overline{\partial}\varphi_{t}\wedge\frac{d}{dt}\overline{\partial}\varphi_{t}
\\
&=&(t\overline{\partial}\varphi''+(1-t)\overline{\varphi}\varphi')\wedge(\overline{\partial}\varphi''-\overline{\partial}\varphi')
\\
&=&t\overline{\partial}\varphi'\wedge\overline{\partial}\varphi''+(1-t)\overline{\partial}\varphi'\wedge\overline{\partial}\varphi''
\ \ = \ \
\overline{\partial}\varphi'\wedge\overline{\partial}\varphi''.
\end{eqnarray*}
If $\varphi'=0$, then
$\overline{\partial}\varphi_{t}\wedge\overline{\partial}\dot{\varphi}_{t}=0$
and hence, by taking the complex conjugate,
$\partial\varphi_{t}\wedge\partial\dot{\varphi}_{t}=0$. This is a
reason why in the Corollary \ref{2.2} there are only three terms.
\end{remark}

Let $S$ be a non-empty set and $A$ an additive group. A mapping
$\mathcal{N}: S\times S\to A$ is said to satisfy the 1-cocycle
condition if
\begin{itemize}

\item[(i)]
$\mathcal{N}(\sigma_{1},\sigma_{2})+\mathcal{N}(\sigma_{2},\sigma_{1})=0$;

\item[(ii)]
$\mathcal{N}(\sigma_{1},\sigma_{2})+\mathcal{N}(\sigma_{2},\sigma_{3})+\mathcal{N}(\sigma_{3},\sigma_{1})=0$.

\end{itemize}

\begin{corollary} \label{cor2.4} (1) The functional $\mathcal{L}^{{\rm M}}_{\omega}$
satisfies the $1$-cocycle condition.\\
(2) For any $\varphi\in\mathcal{P}_{\omega}$ and any constant
$C\in\mathbb{R}$, we have
\begin{equation}
\mathcal{L}^{{\rm
M}}_{\omega}(\varphi,\varphi+C)=C\cdot\left(1-\frac{{\rm
Err}_{\omega}(\varphi)}{V_{\omega}}\right).\label{2.12}
\end{equation}
In particular, if
$\partial\overline{\partial}\omega=\partial\omega\wedge\overline{\partial}\omega=0$,
then $\mathcal{L}^{{\rm M}}_{\omega}(\varphi,\varphi+C)=C$.\\
(3) For any $\varphi_{1},\varphi_{2}\in\mathcal{P}_{\omega}$ and any
constant $C\in\mathbb{R}$, we have
\begin{equation}
\mathcal{L}^{{\rm
M}}_{\omega}(\varphi_{1},\varphi_{2}+C)=\mathcal{L}^{{\rm
M}}_{\omega}(\varphi_{1},\varphi_{2})+C\cdot\left(1-\frac{{\rm
Err}_{\omega}(\varphi_{2})}{V_{\omega}}\right).\label{2.13}
\end{equation}
In particular, if
$\partial\overline{\partial}\omega=\partial\omega\wedge\overline{\partial}\omega=0$,
then $\mathcal{L}^{{\rm
M}}_{\omega}(\varphi_{1},\varphi_{2}+C)=\mathcal{L}^{{\rm
M}}_{\omega}(\varphi_{1},\varphi_{2})+C$.
\end{corollary}

\begin{proof} The proof is similar to that given in \cite{L}.
\end{proof}
\section{Aubin-Yau functionals on compact complex three-folds}

Let $(X,g)$ be a compact complex manifold of the complex dimension
$3$ and $\omega$ be its associated real $(1,1)$-form. We recall some
notation in \cite{L}. For any $\varphi\in\mathcal{P}_{\omega}$ we
set
\begin{eqnarray}
\mathcal{I}^{{\rm AY}}_{\omega|\bullet}(\varphi)&=&
\frac{1}{V_{\omega}}\int_{X}\varphi(\omega^{3}-\omega^{3}_{\varphi}),
\label{3.1}\\
\mathcal{J}^{{\rm
AY}}_{\omega|\bullet}(\varphi)&=&\int^{1}_{0}\frac{\mathcal{I}^{{\rm
AY}}_{\omega|\bullet}(s\cdot\varphi)}{s}ds \ \ = \ \
\frac{1}{V_{\omega}}\int^{1}_{0}\int_{X}\varphi(\omega^{3}-\omega^{3}_{s\cdot\varphi})ds.\label{3.2}
\end{eqnarray}
Two relations showed in \cite{L} are
\begin{eqnarray}
\frac{3}{4}\mathcal{I}^{{\rm
AY}}_{\omega|\bullet}(\varphi)-\mathcal{J}^{{\rm
AY}}_{\omega|\bullet}(\varphi)&=&\frac{1}{V_{\omega}}\int_{X}\varphi
(-\sqrt{-1}\partial\overline{\partial}\varphi)
\wedge\sum^{2}_{j=1}\frac{j}{4}\omega^{2-j}\wedge\omega^{j}_{\varphi},
\label{3.3}\\
4\mathcal{J}^{{\rm AY}}_{\omega|\bullet}(\varphi)-\mathcal{I}^{{\rm
AY}}_{\omega|\bullet}(\varphi)&=&\frac{1}{V_{\omega}}\int_{X}\varphi(-\sqrt{-1}\partial\overline{\partial}\varphi)
\wedge\sum^{1}_{j=0}(2-j)\omega^{2-j}\wedge\omega^{j}_{\varphi}.\label{3.4}
\end{eqnarray}

According to the expression of $\mathcal{L}^{{\rm
M}}_{\omega}(\varphi)$, we set
\begin{eqnarray}
\mathcal{A}_{\omega}(\varphi)&:=&\sum^{1}_{i=0}\frac{i+1}{2V_{\omega}}\int_{X}\varphi
\omega^{i}_{\varphi}\wedge\omega^{1-i}\wedge-\sqrt{-1}\partial\omega\wedge\overline{\partial}\varphi,
\label{3.5}\\
\mathcal{B}_{\omega}(\varphi)&:=&\sum^{1}_{i=0}\frac{i+1}{2V_{\omega}}\int_{X}\varphi
\omega^{i}_{\varphi}\wedge\omega^{1-i}\wedge\sqrt{-1}\overline{\partial}\omega\wedge\partial\varphi.\label{3.6}
\end{eqnarray}

Using (\ref{3.3}) we obtain
\begin{eqnarray*}
& & \frac{3}{4}\mathcal{I}^{{\rm
AY}}_{\omega|\bullet}(\varphi)-\mathcal{J}^{{\rm
AY}}_{\omega|\bullet}(\varphi) \\
&=&\frac{1}{V_{\omega}}\int_{X}\left(\varphi\sum^{2}_{j=1}\frac{j}{4}
\omega^{2-j}\wedge\omega^{j}_{\varphi}\right)\wedge(-\sqrt{-1}\partial\overline{\partial}\varphi)
\\
&=&\frac{1}{V_{\omega}}\int_{X}\sqrt{-1}\partial
\left(\varphi\sum^{2}_{j=1}\frac{j}{4}\omega^{2-j}\wedge\omega^{j}_{\varphi}\right)\wedge\overline{\partial}\varphi
\\
&=&\frac{1}{V_{\omega}}\int_{X}\sqrt{-1}
\left(\partial\varphi\wedge\sum^{2}_{j=1}\frac{j}{4}\omega^{2-j}\wedge\omega^{j}_{\varphi}\right)\wedge\overline{\partial}\varphi
\\
&+&\frac{1}{V_{\omega}}\int_{X}\sqrt{-1}\varphi\sum^{2}_{j=1}\frac{j}{4}[
(2-j)\omega^{1-j}\wedge\partial\omega\wedge\omega^{j}_{\varphi}+\omega^{2-j}\wedge
j\omega^{j-1}_{\varphi}\wedge\partial\omega]\wedge\overline{\partial}\varphi;
\end{eqnarray*}
setting $i=j-1$ in the third term gives
\begin{eqnarray*}
& & \frac{3}{4}\mathcal{I}^{{\rm
AY}}_{\omega|\bullet}(\varphi)-\mathcal{J}^{{\rm
AY}}_{\omega|\bullet}(\varphi) \\
&=&\frac{1}{V_{\omega}}\sum^{2}_{j=1}\frac{j}{4}\int_{X}\sqrt{-1}\partial\varphi\wedge\overline{\partial}\varphi
\wedge\omega^{2-j}\wedge\omega^{j}_{\varphi} \\
&+&\frac{1}{V_{\omega}}\sum^{2}_{j=1}\frac{j(2-j)}{4}\int_{X}\sqrt{-1}\varphi
\omega^{1-j}\wedge\partial\omega\wedge\omega^{j}_{\varphi}\wedge\overline{\partial}\varphi
\\
&+&\frac{1}{V_{\omega}}\sum^{2}_{j=1}\frac{j^{2}}{4}\int_{X}\sqrt{-1}\varphi\omega^{2-j}
\wedge\omega^{j-1}_{\varphi}\wedge\partial\omega\wedge\overline{\partial}\varphi
\\
&=&\frac{1}{V_{\omega}}\sum^{2}_{i=1}\frac{i}{4}\int_{X}\sqrt{-1}\partial\varphi\wedge\overline{\partial}\varphi
\wedge\omega^{2-i}\wedge\omega^{i}_{\varphi}+\frac{1}{4V_{\omega}}\int_{X}\varphi\omega_{\varphi}\wedge
\sqrt{-1}\partial\omega\wedge\overline{\partial}\varphi \\
&+&\frac{1}{V_{\omega}}\sum^{1}_{i=0}\frac{(i+1)^{2}}{4}\int_{X}\varphi\omega^{i}_{\varphi}\wedge\omega^{1-i}\wedge
\sqrt{-1}\partial\omega\wedge\overline{\partial}\varphi \\
&=&\frac{1}{V_{\omega}}\sum^{2}_{i=1}\frac{i}{4}\int_{X}\sqrt{-1}\partial\varphi
\wedge\overline{\partial}\varphi\wedge\omega^{i}_{\varphi}\wedge\omega^{2-i}
\\
&+&\frac{5}{4V_{\omega}}\int_{X}\varphi\omega_{\varphi}
\wedge\sqrt{-1}\partial\omega\wedge\overline{\partial}\varphi
+\frac{1}{4V_{\omega}}\int_{X}\varphi\omega\wedge\sqrt{-1}\partial\omega\wedge\overline{\partial}\varphi.
\end{eqnarray*}
To simplify the notation, we set
\begin{equation}
\mathcal{C}_{\omega}(\varphi):=\frac{3}{4V_{\omega}} \int_{X}\varphi
\omega_{\varphi}\wedge\sqrt{-1}\partial\omega\wedge\overline{\partial}\varphi.\label{3.7}
\end{equation}
Therefore
\begin{eqnarray}
& & \frac{3}{4}\mathcal{I}^{{\rm
AY}}_{\omega|\bullet}(\varphi)-\mathcal{J}^{{\rm
AY}}_{\omega|\bullet}(\varphi) \label{3.8}\\
&=&\frac{1}{V_{\omega}}\sum^{2}_{i=1}\frac{i}{4}\int_{X}\sqrt{-1}\partial\varphi\wedge
\overline{\partial}\varphi\wedge\omega^{i}_{\varphi}\wedge\omega^{2-i}
-\frac{1}{2}\mathcal{A}_{\omega}(\varphi)+\mathcal{C}_{\omega}(\varphi).\nonumber
\end{eqnarray}
On the other hand, using the slightly different method, we obtain
(see \ref{A.1})
\begin{eqnarray}
& & \frac{3}{4}\mathcal{I}^{{\rm
AY}}_{\omega|\bullet}(\varphi)-\mathcal{J}^{{\rm
AY}}_{\omega|\bullet}(\varphi) \label{3.9}\\
&=&\frac{1}{V_{\omega}}\sum^{2}_{i=1}\frac{i}{4}\int_{X}\sqrt{-1}\partial\varphi
\wedge\overline{\partial}\varphi\wedge\omega^{i}_{\varphi}\wedge\omega^{2-i}
-\frac{1}{2}\mathcal{B}_{\omega}(\varphi)+\mathcal{D}_{\omega}(\varphi)\nonumber
\end{eqnarray}
where
\begin{equation}
\mathcal{D}_{\omega}(\varphi)=\frac{3}{4V_{\omega}}
\int_{X}\varphi\omega_{\varphi}\wedge-\sqrt{-1}\overline{\partial}\omega
\wedge\partial\varphi.\label{3.10}
\end{equation}
Equations (\ref{3.8}) and (\ref{3.9}) implies
\begin{eqnarray}
\frac{3}{4}\mathcal{I}^{{\rm
AY}}_{\omega|\bullet}(\varphi)-\mathcal{J}^{{\rm
AY}}_{\omega|\bullet}(\varphi)&=&\frac{1}{V_{\omega}}\sum^{2}_{i=1}\frac{i}{4}\int_{X}\sqrt{-1}\partial\varphi\wedge
\overline{\partial}\varphi\wedge\omega^{i}_{\varphi}\wedge\omega^{2-i}\nonumber
\\
&-&\frac{\mathcal{A}_{\omega}(\varphi)+\mathcal{B}_{\omega}(\varphi)}{4}+\frac{\mathcal{C}_{\omega}(\varphi)+\mathcal{D}_{\omega}(\varphi)}{2}.\label{3.11}
\end{eqnarray}

By the definition we have
\begin{eqnarray*}
\mathcal{J}^{{\rm
AY}}_{\omega|\bullet}(\varphi)&=&\frac{1}{V_{\omega}}\int^{1}_{0}\int_{X}(\varphi\omega^{3}-
\varphi\omega^{3}_{s\varphi})ds \ \ = \ \
\frac{1}{V_{\omega}}\int_{X}\varphi\omega^{3}-\frac{1}{V_{\omega}}\int^{1}_{0}\int_{X}
\varphi\omega^{3}_{t\varphi}dt \\
&=&\frac{1}{V_{\omega}}\int_{X}\varphi\omega^{3}-(\mathcal{L}^{{\rm
M}}_{\omega}(\varphi)-\mathcal{A}_{\omega}(\varphi)-\mathcal{B}_{\omega}(\varphi))
\\
&=&\frac{1}{V_{\omega}}\int_{X}\varphi\omega^{3}-\mathcal{L}^{{\rm
M}}_{\omega}(\varphi)+\mathcal{A}_{\omega}(\varphi)+\mathcal{B}_{\omega}(\varphi).
\end{eqnarray*}

If we define
\begin{eqnarray}
\mathcal{E}_{\omega}(\varphi)&=&\frac{9}{V_{\omega}}
\int_{X}\varphi\omega
\wedge\sqrt{-1}\partial\omega\wedge\overline{\partial}\varphi, \label{3.12}\\
\mathcal{A}^{1}_{\omega}(\varphi)&=&\frac{1}{2V_{\omega}}
\int_{X}\varphi\omega\wedge-\sqrt{-1}\partial
\omega\wedge\overline{\partial}\varphi \label{3.13}\\
\mathcal{A}^{2}_{\omega}(\varphi)&=&\frac{1}{V_{\omega}}
\int_{X}\varphi\omega_{\varphi}\wedge-\sqrt{-1}
\partial\omega\wedge\overline{\partial}\varphi,\label{3.14}
\end{eqnarray}
then
$\mathcal{A}^{1}_{\omega}(\varphi)+\mathcal{A}^{2}_{\omega}(\varphi)=
\mathcal{A}_{\omega}(\varphi)$ and it follows that (see \ref{A.1})
\begin{eqnarray}
4\mathcal{J}^{{\rm AY}}_{\omega|\bullet}-\mathcal{I}^{{\rm
AY}}_{\omega|\bullet}(\varphi)&=&\frac{1}{V_{\omega}}\sum^{1}_{i=0}(2-i)\int_{X}\sqrt{-1}\partial\varphi
\wedge\overline{\partial}\varphi\wedge\omega^{i}_{\varphi}\wedge\omega^{2-i}
\nonumber\\
&+&\mathcal{E}_{\omega}(\varphi)+8\mathcal{A}^{1}_{\omega}(\varphi)-\mathcal{A}^{2}_{\omega}(\varphi).\label{3.15}
\end{eqnarray}

Introduce
\begin{eqnarray}
\mathcal{F}_{\omega}(\varphi)&=&\frac{9}{V_{\omega}}\int_{X} \varphi
\omega\wedge-\sqrt{-1}\overline{\partial}\omega
\wedge\partial\varphi, \label{3.16}\\
\mathcal{B}^{1}_{\omega}(\varphi)&=&\frac{1}{2V_{\omega}}\int_{X}
\varphi\omega\wedge\sqrt{-1}\overline{\partial}\omega\wedge
\partial\varphi \label{3.17}\\
\mathcal{B}^{2}_{\omega}(\varphi)&=&\frac{1}{V_{\omega}}\int_{X}\varphi
\omega_{\varphi}\wedge\sqrt{-1}\overline{\partial}\omega\wedge\partial\varphi.\label{3.18}
\end{eqnarray}
Then
$\mathcal{B}^{1}_{\omega}(\varphi)+\mathcal{B}^{2}_{\omega}(\varphi)=\mathcal{B}_{\omega}(\varphi)$
and hence (see \ref{A.3})
\begin{eqnarray}
4\mathcal{J}^{{\rm AY}}_{\omega|\bullet}(\varphi)-\mathcal{I}^{{\rm
AY}}_{\omega|\bullet}(\varphi)&=&\frac{1}{V_{\omega}}\sum^{1}_{i=0}(2-i)
\int_{X}\sqrt{-1}\partial\varphi\wedge\overline{\partial}\varphi\wedge\omega^{i}_{\varphi}
\wedge\omega^{2-i} \nonumber \\
&+&\mathcal{F}_{\omega}(\varphi)+8\mathcal{B}^{1}_{\omega}(\varphi)-\mathcal{B}^{2}_{\omega}(\varphi).\label{3.19}
\end{eqnarray}
(\ref{3.15}) and (\ref{3.19}) together gives
\begin{eqnarray}
4\mathcal{J}^{{\rm AY}}_{\omega|\bullet}(\varphi)-\mathcal{I}^{{\rm
AY}}_{\omega|\bullet}(\varphi)&=&\frac{1}{V_{\omega}}\sum^{1}_{i=0}(2-i)
\int_{X}\sqrt{-1}\partial\varphi\wedge\overline{\partial}\varphi
\wedge\omega^{i}_{\varphi}\wedge\omega^{n-1-i}\nonumber \\
&+&\frac{\mathcal{E}_{\omega}(\varphi)+\mathcal{F}_{\omega}(\varphi)}{2}
+4[\mathcal{A}^{1}_{\omega}(\varphi)+\mathcal{B}^{1}_{\omega}(\varphi)]\label{3.20}\\
&-&\frac{\mathcal{A}^{2}_{\omega}(\varphi)
+\mathcal{B}^{2}_{\omega}(\varphi)}{2}.\nonumber
\end{eqnarray}

Now, we define Aubin-Yau functionals over any compact complex
manifolds as follows:
\begin{eqnarray}
\mathcal{I}^{{\rm AY}}_{\omega}(\varphi)&:=&\mathcal{I}^{{\rm
AY}}_{\omega|\bullet}(\varphi)\nonumber\\
&+&a^{1}_{1}\mathcal{A}^{1}_{\omega}(\varphi)+a^{2}_{1}\mathcal{A}^{2}_{\omega}(\varphi)
+b^{1}_{1}\mathcal{B}^{1}_{\omega}(\varphi)+b^{2}_{1}\mathcal{B}^{2}_{\omega}(\varphi)\nonumber\\
&+&c_{1}\mathcal{C}_{\omega}(\varphi)+d_{1}\mathcal{D}_{\omega}(\varphi)
+e_{1}\mathcal{E}_{\omega}(\varphi)+f_{1}\mathcal{F}_{\omega}(\varphi),\label{3.21}
\\
\mathcal{J}^{{\rm AY}}_{\omega}(\varphi)&:=&-\mathcal{L}^{{\rm
M}}_{\omega}(\varphi)+\frac{1}{V_{\omega}}\int_{X}\varphi\omega^{3}\nonumber\\
&+&a^{1}_{2}\mathcal{A}^{1}_{\omega}(\varphi)+a^{2}_{2}\mathcal{A}^{2}_{\omega}(\varphi)
+b^{1}_{2}\mathcal{B}^{1}_{\omega}(\varphi)+b^{2}_{2}\mathcal{B}^{2}_{\omega}(\varphi)\nonumber\\
&+&c_{2}\mathcal{C}_{\omega}(\varphi)+d_{2}\mathcal{D}_{\omega}(\varphi)
+e_{2}\mathcal{E}_{\omega}(\varphi)+f_{2}\mathcal{F}_{\omega}(\varphi),
\nonumber\\
&=&\mathcal{J}^{{\rm AY}}_{\omega|\bullet}(\varphi) \nonumber \\
&+&(a^{1}_{2}-1)\mathcal{A}^{1}_{\omega}(\varphi)+(a^{2}_{2}-1)\mathcal{A}^{2}_{\omega}(\varphi)
+(b^{1}_{2}-1)\mathcal{B}^{1}_{\omega}(\varphi)+(b^{2}_{2}-1)\mathcal{B}_{\omega}^{2}(\varphi)
\nonumber\\
&+&c_{2}\mathcal{C}_{\omega}(\varphi)+d_{2}\mathcal{D}_{\omega}(\varphi)
+e_{2}\mathcal{E}_{\omega}(\varphi)+f_{2}\mathcal{F}_{\omega}(\varphi).\label{3.22}
\end{eqnarray}

Plugging (\ref{3.21}) and (\ref{3.22}) into (\ref{3.20}) and
(\ref{3.11}), we obtain
\begin{equation}
\frac{3}{4}\mathcal{I}^{{\rm
AY}}_{\omega}(\varphi)-\mathcal{J}^{{\rm
AY}}_{\omega}(\varphi)=\frac{1}{V_{\omega}}\sum^{2}_{i=1}\frac{i}{4}\int_{X}\sqrt{-1}\partial\varphi
\wedge\overline{\partial}\varphi\wedge\omega^{i}_{\varphi}\wedge\omega^{2-i}\geq0\label{3.23}
\end{equation}
and
\begin{equation}
4\mathcal{J}^{{\rm AY}}_{\omega}(\varphi)-\mathcal{I}^{{\rm
AY}}_{\omega}(\varphi)=\sum^{2}_{i=0}\frac{2-i}{V_{\omega}}\int_{X}\sqrt{-1}\partial
\varphi\wedge\overline{\partial}\varphi\wedge\omega^{i}_{\varphi}\wedge\omega^{2-i}\geq0\label{3.24}
\end{equation}
where we require that constants satisfy the following linear
equations system
\begin{eqnarray}
\frac{3}{4}a^{1}_{1}-(a^{1}_{2}-1)&=&\frac{1}{4},\label{3.25} \ \ \
\frac{3}{4}a^{2}_{1}-(a^{2}_{2}-1) \ \ = \ \ \frac{1}{4}, \\
\frac{3}{4}b^{1}_{1}-(b^{1}_{2}-1)&=&\frac{1}{4}, \label{3.26} \ \ \
\frac{3}{4}b^{2}_{1}-(b^{2}_{2}-1) \ \ = \ \ \frac{1}{4},
\\
\frac{3}{4}c_{1}-c_{2}&=&-\frac{1}{2}, \label{3.27} \ \ \
\frac{3}{4}d_{1}-d_{2} \ \ = \ \ -\frac{1}{2}, \\
\frac{3}{4}e_{1}-e_{2}&=&0,\label{3.28} \ \ \
\frac{3}{4}f_{1}-f_{2} \ \ = \ \ 0, \\
4(a^{1}_{2}-1)-a^{1}_{1}&=&-4,\label{3.29} \ \ \
4(a^{2}_{2}-1)-a^{2}_{1} \ \ = \ \ \frac{1}{2}, \\
4(b^{1}_{2}-1)-b^{1}_{1}&=&-4, \label{3.30} \ \ \
4(b^{2}_{2}-1)-b^{2}_{1} \ \ = \ \ \frac{1}{2}, \\
4c_{2}-c_{1}&=&0, \label{3.31} \ \ \
4d_{2}-d_{1} \ \ = \ \ 0, \\
4e_{2}-e_{1}&=&-\frac{1}{2},\label{3.32} \ \ \ 4f_{2}-f_{1} \ \ = \
\ -\frac{1}{2}.
\end{eqnarray}
The constants $a^{j}_{i},b^{j}_{i},c_{i},d_{i},e_{i}$ and $f_{i}$ ,
calculated in Appendix B, are
\begin{eqnarray}
a^{1}_{1}&=&b^{1}_{1} \ \ = \ \ -\frac{3}{2}, \ \ \ a^{1}_{2} \ \ =
\ \
b^{1}_{2} \ \ = \ \ -\frac{3}{8},\label{3.33}\\
a^{2}_{1}&=&b^{2}_{1} \ \ = \ \ \frac{3}{4}, \ \ \ a^{2}_{2} \ \ = \
\ b^{2}_{2} \ \ = \ \
\frac{3}{4}\left(1+\frac{3}{4}\right) \ \ = \ \ \frac{21}{16},\label{3.34}\\
c_{1}&=&d_{1} \ \ = \ \ -1, \ \ \ e_{2} \ \ = \ \
f_{2} \ \ = \ \ -\frac{3}{16}, \label{3.35}\\
c_{2}&=&d_{2} \ \ = \ \ e_{1} \ \ = \ \ f_{1} \ \ = \ \
-\frac{1}{4}.\label{3.36}
\end{eqnarray}

The explicit formulas for $\mathcal{I}^{{\rm AY}}_{\omega}(\varphi)$
and $\mathcal{J}^{{\rm AY}}_{\omega}(\varphi)$ are given in
Proposition \ref{C.1} and \ref{C.2} respectively. Namely,
\begin{eqnarray}
& & \mathcal{I}^{{\rm AY}}_{\omega}(\varphi) \ \ = \ \
\frac{1}{V_{\omega}}\int_{X}\varphi(\omega^{3}-\omega^{3}_{\varphi})\label{3.37}\\
&-&\frac{3}{2V_{\omega}}\int_{X}\varphi\omega_{\varphi}
\wedge\sqrt{-1}\partial\omega\wedge\overline{\partial}\varphi-
\frac{3}{2V_{\omega}}\int_{X}\varphi\omega\wedge\sqrt{-1}\partial\omega\wedge
\overline{\partial}\varphi \nonumber\\
&+&\frac{3}{2V_{\omega}}\int_{X}\varphi\omega_{\varphi}
\wedge\sqrt{-1}\overline{\partial}\omega\wedge\partial\varphi
+\frac{3}{2V_{\omega}}\int_{X}\varphi\omega\wedge\sqrt{-1}\overline{\partial}\omega\wedge\partial\varphi
\nonumber,\\
& & \mathcal{J}^{{\rm AY}}_{\omega}(\varphi) \ \ = \ \
-\mathcal{L}^{{\rm
M}}_{\omega}(\varphi)+\frac{1}{V_{\omega}}\int_{X}\varphi\omega^{3}
\label{3.38}\\
&-&\frac{3}{2V_{\omega}}\int_{X}\varphi\omega_{\varphi}
\wedge\sqrt{-1}\partial\omega\wedge\overline{\partial}\varphi
-\frac{3}{2V_{\omega}}\int_{X}\varphi\omega\wedge\sqrt{-1}\partial\omega\wedge\overline{\partial}\varphi\nonumber\\
&+&\frac{3}{2V_{\omega}}\int_{X}\varphi\omega_{\varphi}
\wedge\sqrt{-1}\overline{\partial}\omega\wedge\partial\varphi
+\frac{3}{2V_{\omega}}\int_{X}\varphi\omega\wedge\sqrt{-1}\overline{\partial}\omega\wedge\partial\varphi.\nonumber
\end{eqnarray}

From (\ref{3.23}), (\ref{3.24}), (\ref{3.31}) and (\ref{3.38}), we
deduce the following

\begin{theorem} \label{thm3.1} For any $\varphi\in\mathcal{P}_{\omega}$, one has
\begin{eqnarray}
\frac{3}{4}\mathcal{I}^{{\rm
AY}}_{\omega}(\varphi)-\mathcal{J}^{{\rm
AY}}_{\omega}(\varphi)&\geq&0, \label{3.39}\\
4\mathcal{J}^{{\rm AY}}_{\omega}(\varphi)-\mathcal{I}^{{\rm
AY}}_{\omega}(\varphi)&\geq&0.\label{3.40}
\end{eqnarray}
In particular
\begin{eqnarray*}
\frac{1}{4}\mathcal{I}^{{\rm
AY}}_{\omega}(\varphi)&\leq&\mathcal{J}^{{\rm AY}}_{\omega}(\varphi)
\ \ \leq \ \ \frac{3}{4}\mathcal{I}^{{\rm AY}}_{\omega}(\varphi),
\\
\frac{3}{4}\mathcal{J}^{{\rm
AY}}_{\omega}(\varphi)&\leq&\mathcal{I}^{{\rm AY}}_{\omega}(\varphi)
\ \ \leq \ \ 4\mathcal{J}^{{\rm AY}}_{\omega}(\varphi), \\
\frac{1}{3}\mathcal{J}^{{\rm AY}}_{\omega}(\varphi) \ \ \leq \ \
\frac{1}{4}\mathcal{J}^{{\rm AY}}_{\omega}(\varphi)&\leq&
\mathcal{I}^{{\rm
AY}}_{\omega}(\varphi)-\mathcal{J}^{{\rm AY}}_{\omega}(\varphi)\\
&\leq&\frac{3}{4}\mathcal{I}^{{\rm AY}}_{\omega}(\varphi) \ \ \leq \
\ n\mathcal{J}^{{\rm AY}}_{\omega}(\varphi).
\end{eqnarray*}
\end{theorem}

\section{Volume estimates}
Let $(X,g)$ be a compact Hermitian manifold of the complex dimension
$n$ and $\omega$ be its associated real $(1,1)$-form. Define
\begin{equation}
\mathcal{P}_{\omega}:=\{\varphi\in
C^{\infty}(X)_{\mathbb{R}}|\omega_{\varphi}:=\omega+\sqrt{-1}\partial\overline{\partial}\varphi>0\},\label{4.1}
\end{equation}
and
$\mathcal{P}^{0}_{\omega}:=\left\{\varphi\in\mathcal{P}_{\omega}\Big|\sup_{X}\varphi=0\right\}$.
Consider the quantity
\begin{equation}
{\rm InfErr}_{\omega}:=\inf_{\varphi\in\mathcal{P}^{0}_{\omega}}{\rm
Err}_{\omega}(\varphi),\label{4.2}
\end{equation}
where
\begin{equation}
{\rm
Err}_{\omega}(\varphi):=\int_{X}\omega^{n}-\int_{X}\omega^{n}_{\varphi}.\label{4.3}
\end{equation}
It's clear that ${\rm InfErr}_{\omega}\leq0$. But we don't know
whether the quantity ${\rm InfErr}_{\omega}$ is finite. For $n=2$,
V. Tosatti and B. Weinkove \cite{TW1} showed that ${\rm
InfErr}_{\omega}$ is always bounded from below, using the existence
of Ganduchon metrics on any compact Hermitian manifolds. If $\omega$
satisfies the condition (see \cite{GL}, \cite{TW1})
\begin{equation}
\partial\overline{\partial}(\omega^{k})=0, \ \ \ k=1,2,\label{4.4}
\end{equation}
one can show that the quantity ${\rm InfErr}_{\omega}=0$ (see
\cite{DK}, \cite{L}, or \cite{TW1}).

When $n=2$, the condition (\ref{4.4}) reduces to
\begin{equation}
\partial\overline{\partial}\omega=0,\label{4.5}
\end{equation}
which is a Ganduchon metric; however, if $n=3$, we can show that
${\rm InfErr}_{\omega}$ is bounded from below under this condition.

\begin{theorem} \label{thm4.1} Suppose that $(X,g)$ is a compact Hermitian manifold
of the complex dimension $3$ and $\omega$ is its associated real
$(1,1)$-form. If $\partial\overline{\partial}\omega=0$, then ${\rm
InfErr}_{\omega}$ is bounded from below. More precisely, we have
\begin{equation}
{\rm InfErr}_{\omega}\geq3\left(1-e^{2\cdot{\rm
osc}(u)}\right)\cdot\int_{X}\omega^{3}.\label{4.6}
\end{equation}
Here $u$ is a real-valued smooth function on $X$ such that
$\omega_{G}=e^{u}\cdot \omega$ is a Gauduchon metric, i.e.,
$\partial\overline{\partial}(\omega^{2}_{G})=0$.
\end{theorem}

\begin{proof} As in \cite{TW1}, page 21, we compute
\begin{eqnarray*}
\int_{X}\omega^{3}_{\varphi}&=&\int_{X}(\omega^{3}+3\omega^{2}\wedge\sqrt{-1}\partial\overline{\partial}\varphi+3\omega\wedge(\sqrt{-1}\partial\overline{\partial}\varphi)^{2})
\\
&=&\int_{X}(-2\omega^{3}+3\omega^{2}\wedge(\omega+\sqrt{-1}\partial\overline{\partial}\varphi)+3\omega\wedge(\sqrt{-1}\partial\overline{\partial}\varphi)^{2}).
\end{eqnarray*}
Since $\partial\overline{\partial}\omega=0$, the last integral
vanishes, and hence
\begin{eqnarray*}
\int_{X}\omega^{3}_{\varphi}&\leq&\int_{X}-2\omega^{3}+\int_{X}3\left(e^{u-\inf_{X}(u)}\omega\right)^{2}\wedge(\omega+\sqrt{-1}\partial\overline{\partial}\varphi)
\\
&=&\int_{X}-2\omega^{3}+3e^{2(\sup_{X}(u)-\inf_{X}(u))}\int_{X}\omega^{3}\\
&+&\int_{X}3e^{-2\inf_{X}(u)}\cdot\omega_{G}^{2}\wedge\sqrt{-1}\partial\overline{\partial}\varphi
\\
&=&\int_{X}-2\omega^{3}+3e^{2\cdot{\rm osc}(u)}\int_{X}\omega^{3} \
\ = \ \ \left(3e^{2\cdot{\rm osc}(u)}-2\right)\int_{X}\omega^{3}.
\end{eqnarray*}
From the definition of ${\rm InfErr}_{\omega}$, we immediately
obtain
\begin{equation*}
{\rm InfErr}_{\omega}\geq\int_{X}\omega^{3}-\left(3e^{2\cdot{\rm
osc}(u)}-2\right)\int_{X}\omega^{3}=3\left(1-e^{2\cdot{\rm
osc}(u)}\right)\cdot\int_{X}\omega^{3}
\end{equation*}
where ${\rm osc}(u):=\sup_{X}(u)-\inf_{X}(u)$.
\end{proof}

\begin{theorem} \label{thm4.2} Suppose that $(X,g)$ is a compact Hermitian manifold
of the complex dimension $3$ and $\omega$ is its associated real
$(1,1)$-form. We select a real-valued smooth function $u$ on $X$ so
that $e^{u}\cdot \omega$ is a Gauduchon metric. If
\begin{equation*}
{\rm osc}(u)=\sup_{X}(u)-\inf_{X}(u)\leq\frac{1}{2}\cdot{\rm
ln}\frac{3}{2}, \ \ \
\partial\overline{\partial}\omega=0,
\end{equation*}
then
\begin{equation}
\inf_{\varphi\in\mathcal{P}^{0}_{\omega}}\int_{X}\omega^{3}_{\varphi}\geq\int_{X}\omega^{3}>0.\label{4.7}
\end{equation}
\end{theorem}

\begin{proof} Using the similar procedure, we deduce
\begin{equation*}
\int_{X}\omega^{3}_{\varphi}\geq\left(3\cdot
e^{2(\inf_{X}(u)-\sup_{X}(u))}-2\right)\int_{X}\omega^{3}.
\end{equation*}
Since $\sup_{X}u-\inf_{X}u\leq\frac{1}{2}\cdot{\rm ln}\frac{3}{2}$,
it follows that $3\cdot e^{2(\inf_{X}u-\sup_{X}u)}-2\geq1$.
\end{proof}

\appendix
\section{Proof the identities (\ref{3.8}), (\ref{3.15}) and (\ref{3.19})}
In Appendix A we verify the identities (\ref{3.8}), (\ref{3.15}) and
(\ref{3.19}).

\begin{eqnarray}
& &\frac{3}{4}\mathcal{I}^{{\rm
AY}}_{\omega|\bullet}(\varphi)-\mathcal{J}^{{\rm
AY}}_{\omega|\bullet}(\varphi)
\label{A.1}\\
&=&\frac{1}{V_{\omega}}\int_{X}-\sqrt{-1}\overline{\partial}\left(\varphi\sum^{2}_{j=1}
\frac{j}{4}\omega^{2-j}\wedge\omega^{j}_{\varphi}\right)\wedge\partial\varphi
\nonumber\\
&=&\frac{1}{V_{\omega}}\int_{X}-\sqrt{-1}\left(\overline{\partial}\varphi\wedge
\sum^{2}_{j=1}\frac{j}{4}\omega^{2-j}\wedge\omega^{j}_{\varphi}\right)\wedge\partial\varphi
\nonumber\\
&+&\frac{1}{V_{\omega}}\int_{X}-\sqrt{-1}\varphi\sum^{2}_{j=1}\frac{j}{4}[(2-j)\omega^{1-j}
\wedge\overline{\partial}\omega\wedge\omega^{j}_{\varphi}+\omega^{2-j}\wedge
j\omega^{j-1}_{\varphi}\wedge\overline{\partial}\omega]\wedge\partial\varphi
\nonumber\\
&=&\frac{1}{V_{\omega}}\sum^{2}_{i=1}\frac{i}{4}\int_{X}\sqrt{-1}\partial\varphi
\wedge\overline{\partial}\varphi\wedge\omega^{i}_{\varphi}\wedge\omega^{2-i}
\nonumber\\
&+&\frac{1}{V_{\omega}}\sum^{2}_{j=1}\frac{j(2-j)}{4}\int_{X}-\sqrt{-1}\varphi
\omega^{1-j}\wedge\omega^{j}_{\varphi}\wedge\overline{\partial}\omega\wedge\partial\varphi
\nonumber\\
&+&\frac{1}{V_{\omega}}\sum^{2}_{j=1}\frac{j^{2}}{4}\int_{X}-\sqrt{-1}\varphi
\omega^{2-j}\wedge\omega^{j-1}_{\varphi}\wedge\overline{\partial}\omega\wedge\partial\varphi
\nonumber\\
&=&\frac{1}{V_{\omega}}\sum^{2}_{i=1}\frac{i}{4}\int_{X}\sqrt{-1}\partial\varphi\wedge\overline{\partial}
\varphi\wedge\omega^{2-i}\wedge\omega^{i}_{\varphi}+\frac{1}{4V_{\omega}}\int_{X}-\sqrt{-1}\varphi
\omega_{\varphi}\wedge\overline{\partial}\omega\wedge\partial\varphi
\nonumber\\
&+&\frac{1}{V_{\omega}}\sum^{1}_{i=0}\frac{(i+1)^{2}}{4}\int_{X}-\sqrt{-1}\varphi
\omega^{1-i}\wedge\omega^{i}_{\varphi}\wedge\overline{\partial}\omega\wedge\partial\varphi
\nonumber\\
&=&\frac{1}{V_{\omega}}\sum^{2}_{i=1}\frac{i}{4}\int_{X}\sqrt{-1}\partial\varphi
\wedge\overline{\partial}\varphi\wedge\omega^{2-i}\wedge\omega^{i}_{\varphi}
-\frac{1}{2}\mathcal{B}_{\omega}(\varphi)+\mathcal{D}_{\omega}(\varphi)\nonumber
\end{eqnarray}
which gives (\ref{3.9}). Calculate
\begin{eqnarray*}
& & 4\mathcal{J}^{{\rm AY}}_{\omega|\bullet}-\mathcal{I}^{{\rm
AY}}_{\omega|\bullet}(\varphi) \label{A.2}\\
&=&\frac{1}{V_{\omega}}\int_{X}\sqrt{-1}\partial\left(\varphi\sum^{2}_{j=0}(2-j)\omega^{2-j}\wedge\omega^{j}_{\varphi}\right)
\wedge\overline{\partial}\varphi \\
&=&\frac{1}{V_{\omega}}\int_{X}\sqrt{-1}\partial\varphi\wedge\sum^{2}_{j=0}
(2-j)\omega^{2-j}\wedge\omega^{j}_{\varphi}\wedge\overline{\partial}\varphi
\\
&+&\frac{1}{V_{\omega}}\int_{X}\sqrt{-1}\varphi\sum^{2}_{j=0}[(2-j)^{2}\omega^{1-j}
\wedge\partial\omega\wedge\omega^{j}_{\varphi}+(2-j)j\omega^{2-j}\wedge\omega^{j-1}_{\varphi}\wedge\partial\omega]\wedge\overline{\partial}\varphi
\\
&=&\frac{1}{V_{\omega}}\sum^{2}_{i=0}(2-i)\int_{X}\sqrt{-1}\partial\varphi
\wedge\overline{\partial}\varphi\wedge\omega^{i}_{\varphi}\wedge\omega^{2-i}
\\
&+&\frac{1}{V_{\omega}}\sum^{1}_{j=0}(2-j)^{2}\int_{X}\varphi\omega^{1-j}\wedge\omega^{j}_{\varphi}\wedge\sqrt{-1}\partial\omega
\wedge\overline{\partial}\varphi+\frac{1}{V_{\omega}}\int_{X}\varphi\omega\wedge\sqrt{-1}
\partial\omega\wedge\overline{\partial}\varphi \\
&=&\frac{1}{V_{\omega}}\sum^{2}_{i=0}(2-i)\int_{X}\sqrt{-1}\partial\varphi
\wedge\overline{\partial}\varphi\wedge\omega^{i}_{\varphi}\wedge\omega^{2-i}
\\
&+&\frac{5}{V_{\omega}}\int_{X}\varphi\omega
\wedge\sqrt{-1}\partial\omega\wedge\overline{\partial}\varphi+\frac{1}{V_{\omega}}\int_{X}\varphi\omega_{\varphi}\wedge\sqrt{-1}\partial\omega
\wedge\overline{\partial}\varphi.
\end{eqnarray*}
Using the definitions of $\mathcal{E}_{\omega}(\varphi),
\mathcal{A}^{1}_{\omega}(\varphi)$,
$\mathcal{A}^{2}_{\omega}(\varphi)$, we have
$\mathcal{A}^{1}_{\omega}(\varphi)+\mathcal{A}^{2}_{\omega}(\varphi)=
\mathcal{A}_{\omega}(\varphi)$ and hence (\ref{3.15}) holds.
Similarly, we have
$\mathcal{B}^{1}_{\omega}(\varphi)+\mathcal{B}^{2}_{\omega}(\varphi)=\mathcal{B}_{\omega}(\varphi)$
and
\begin{eqnarray*}
& & 4\mathcal{J}^{{\rm
AY}}_{\omega|\bullet}(\varphi)-\mathcal{I}^{{\rm
AY}}_{\omega|\bullet}(\varphi) \\
&=&\frac{1}{V_{\omega}}\int_{X}-\sqrt{-1}\overline{\partial}\left(\varphi
\sum^{2}_{j=0}(2-j)\omega^{2-j}\wedge\omega^{j}_{\varphi}\right)\wedge\partial\varphi
\\
&=&\frac{1}{V_{\omega}}\int_{X}-\sqrt{-1}\overline{\partial}\varphi\wedge
\sum^{2}_{j=0}(2-j)\omega^{2-j}\wedge\omega^{j}_{\varphi}\wedge\partial\varphi
\\
&+&\frac{1}{V_{\omega}}\int_{X}-\sqrt{-1}\varphi\sum^{2}_{j=0}(2-j)
[(2-j)\omega^{1-j}\wedge\overline{\partial}\omega\wedge\omega^{j}_{\varphi}+j\omega^{2-j}\wedge\omega^{j-1}_{\varphi}\wedge\overline{\partial}\omega]
\wedge\partial\varphi \\
&=&\frac{1}{V_{\omega}}\sum^{2}_{i=0}(2-i)\int_{X}\sqrt{-1}
\partial\varphi\wedge\overline{\partial}\varphi\wedge\omega^{2-i}\wedge
\omega^{i}_{\varphi} \\
&+&\frac{1}{V_{\omega}}\sum^{1}_{i=0}(2-i)^{2}\int_{X}\varphi\omega^{1-i}
\wedge\omega^{i}_{\varphi}\wedge(-\sqrt{-1}\overline{\partial}\omega\wedge\partial\varphi)
+\frac{1}{V_{\omega}}\int_{X}\varphi\omega\wedge(-\sqrt{-1}\overline{\partial}\omega\wedge\partial\varphi)
\\
&=&\frac{1}{V_{\omega}}\sum^{2}_{i=0}(2-i)\int_{X}\sqrt{-1}\partial\varphi\wedge\overline{\partial}\varphi
\wedge\omega^{2-i}\wedge\omega^{i}_{\varphi} \\
&+&\frac{5}{V_{\omega}}\int_{X}\varphi
\omega\wedge(-\sqrt{-1}\overline{\partial}\omega\wedge\partial\varphi)
+\frac{1}{V_{\omega}}\int_{X}\varphi\omega_{\varphi}\wedge(-\sqrt{-1}
\overline{\partial}\omega\wedge\partial\varphi).
\end{eqnarray*}
and hence
\begin{eqnarray}
4\mathcal{J}^{{\rm AY}}_{\omega|\bullet}(\varphi)-\mathcal{I}^{{\rm
AY}}_{\omega|\bullet}(\varphi)&=&\frac{1}{V_{\omega}}\sum^{2}_{i=0}(2-i)
\int_{X}\sqrt{-1}\partial\varphi\wedge\overline{\partial}\varphi\wedge\omega^{i}_{\varphi}
\wedge\omega^{2-i} \nonumber \\
&+&\mathcal{F}_{\omega}(\varphi)+
8\mathcal{B}^{1}_{\omega}(\varphi)-\mathcal{B}^{2}_{\omega}(\varphi).\label{A.3}
\end{eqnarray}
Therefore (\ref{3.15}) and (\ref{3.19}) together gives
\begin{eqnarray}
4\mathcal{J}^{{\rm AY}}_{\omega|\bullet}(\varphi)-\mathcal{I}^{{\rm
AY}}_{\omega|\bullet}(\varphi)&=&\frac{1}{V_{\omega}}\sum^{2}_{i=0}(2-i)
\int_{X}\sqrt{-1}\partial\varphi\wedge\overline{\partial}\varphi
\wedge\omega^{i}_{\varphi}\wedge\omega^{2-i}\nonumber \\
&+&\frac{\mathcal{E}_{\omega}(\varphi)+\mathcal{F}_{\omega}(\varphi)}{2}
+4(\mathcal{A}^{1}_{\omega}(\varphi)+\mathcal{B}^{1}_{\omega}(\varphi))-\frac{\mathcal{A}^{2}_{\omega}(\varphi)
+\mathcal{B}^{2}_{\omega}(\varphi)}{2}.\label{A.3}
\end{eqnarray}

\section{Solve the linear equations system} \label{oldnew}
In this section we try to solve the linear equations system
(\ref{3.25})-(\ref{3.32}). Firstly we solve (\ref{3.25}) and
(\ref{3.29}) as follows: (\ref{3.25}) and (\ref{3.29}) gives us the
following equations
\begin{eqnarray}
\frac{3}{4}a^{1}_{1}-\frac{1}{4}&=&a^{1}_{2}-1, \label{B.1} \ \ \
4(a^{1}_{2}-1)+4 \ \ = \ \ a^{1}_{1}, \\
\frac{3}{4}a^{2}_{1}-\frac{1}{4}&=&a^{2}_{2}-1,\label{B.2} \ \ \
4(a^{2}_{2}-1)-\frac{1}{2} \ \ = \ \ a^{2}_{1}.
\end{eqnarray}
Plugging the first equation into second equation in (\ref{B.1}), we
have
\begin{equation*}
4\left(\frac{3}{4}a^{1}_{1}-\frac{1}{4}\right)+4=a^{1}_{1}
\end{equation*}
which implies
\begin{equation}
a^{1}_{1}=-\frac{3}{2}, \ \ \ a^{1}_{2}=-\frac{3}{8}.\label{B.3}
\end{equation}
Similarly,
\begin{equation*}
4\left(\frac{3}{4}a^{2}_{1}-\frac{1}{4}\right)-\frac{1}{2}=a^{2}_{1},
\end{equation*}
therefore
\begin{equation}
a^{2}_{1}=\frac{3}{4}, \ \ \
a^{2}_{2}=\frac{3}{4}\left(1+\frac{3}{4}\right)=\frac{21}{16}.\label{B.4}
\end{equation}
Secondly, (\ref{3.26}) and (\ref{3.30}) implies
\begin{eqnarray}
\frac{3}{4}b^{1}_{1}-\frac{1}{4}&=&b^{1}_{2}-1,\label{B.5} \ \ \
4(b^{1}_{2}-1) \ \ = \ \ b^{1}_{1}-4, \\
\frac{3}{4}b^{2}_{1}-\frac{1}{4}&=&b^{2}_{2}-1, \label{B.6} \ \ \
4(b^{2}_{2}-1) \ \ = \ \ b^{2}_{1}+\frac{1}{2}.
\end{eqnarray}
The above linear equations system gives
\begin{equation*}
4\left(\frac{3}{4}b^{1}_{1}-\frac{1}{4}\right)=b^{1}_{1}-4, \ \ \
4\left(\frac{3}{4}b^{2}_{1}-\frac{1}{4}\right)=b^{2}_{1}+\frac{1}{2},
\end{equation*}
respectively. Hence
\begin{eqnarray}
b^{1}_{1}&=&-\frac{3}{2}, \label{B.7} \ \ \ b^{1}_{2} \ \ = \ \
-\frac{3}{8}, \\
b^{2}_{1}&=&\frac{3}{4}, \label{B.8} \ \ \
 b^{2}_{2} \ \ = \ \
\frac{3}{4}\left(1+\frac{3}{4}\right) \ \ = \ \ \frac{21}{16}.
\end{eqnarray}
Continuously, equations (\ref{3.27}) and (\ref{3.31}) shows that
\begin{eqnarray*}
\frac{3}{4}c_{1}-c_{2}&=&-\frac{1}{2}, \ \ \ 4c_{2}-c_{1} \ \
= \ \ 0, \\
\frac{3}{4}d_{1}-d_{2}&=&-\frac{1}{2}, \ \ \ 4d_{2}-d_{1} \ \ = \ \
0.
\end{eqnarray*}
Eliminating $c_{2}$ and $d_{2}$ respectively, we have
\begin{equation*}
4\left(\frac{3}{4}c_{1}+\frac{1}{2}\right)-c_{1}=0, \ \ \
4\left(\frac{3}{4}d_{1}+\frac{1}{2}\right)-d_{1}=0.
\end{equation*}
Thus
\begin{eqnarray}
c_{1}&=&-1, \label{B.9} \ \ \ c_{2} \ \ = \ \ -\frac{1}{4} \\
d_{1}&=&-1, \label{B.10} \ \ \ d_{2} \ \ = \ \ -\frac{1}{4}.
\end{eqnarray}
Similarly, from (\ref{3.28}) and (\ref{3.32}) we obtain
\begin{eqnarray*}
\frac{3}{4}e_{1}-e_{2}=0, \ \ \ 4e_{2}-e_{1} \ \ = \ \
-\frac{1}{2}, \\
\frac{3}{4}f_{1}-f_{2}=0, \ \ \ 4f_{2}-f_{1} \ \ = \ \ -\frac{1}{2},
\end{eqnarray*}
and hence
\begin{eqnarray}
e_{1}&=&f_{1} \ \ = \ \ -\frac{1}{4},\label{B.11}\\
e_{2}&=&f_{2} \ \ = \ \ -\frac{3}{16}.\label{B.12}
\end{eqnarray}
\section{Explicit formulas of $\mathcal{I}^{{\rm
AY}}_{\omega}(\varphi)$ and $\mathcal{J}^{{\rm
AY}}_{\omega}(\varphi)$} \label{oldnew}
In this section we give the explicit formulas of $\mathcal{I}^{{\rm
AY}}_{\omega}(\varphi)$ and $\mathcal{J}^{{\rm
AY}}_{\omega}(\varphi)$. Using the constants determined in Appendix
B, we have
\begin{eqnarray*}
\mathcal{I}^{{\rm
AY}}_{\omega}(\varphi)&=&\frac{1}{V_{\omega}}\int_{X}\varphi(\omega^{3}-\omega^{3}_{\varphi})
\\
&+&\frac{3}{4V_{\omega}}\int_{X}\varphi\omega\wedge\sqrt{-1}\partial\omega\wedge\overline{\partial}\varphi
-\frac{3}{4V_{\omega}}\int_{X}\varphi\omega_{\varphi}\wedge\sqrt{-1}\partial\omega\wedge\overline{\partial}\varphi
\\
&-&\frac{3}{4V_{\omega}}\int_{X}\varphi\omega\wedge\sqrt{-1}\overline{\partial}\omega\wedge\partial\varphi
+\frac{3}{4V_{\omega}}\int_{X}\varphi\omega_{\varphi}\wedge\sqrt{-1}\overline{\partial}\omega\wedge\partial\varphi
\\
&-&\frac{3}{4V_{\omega}}\int_{X}\varphi\omega_{\varphi}\wedge\sqrt{-1}\partial\omega
\wedge\overline{\partial}\varphi+\frac{3}{4V_{\omega}}\int_{X}\varphi\omega_{\varphi}\wedge\sqrt{-1}\overline{\partial}\omega\wedge\partial\varphi
\\
&-&\frac{9}{4V_{\omega}}\int_{X}\varphi\omega\wedge\sqrt{-1}\partial\omega\wedge
\overline{\partial}\varphi+
\frac{9}{4V_{\omega}}\int_{X}\varphi\omega\wedge\sqrt{-1}\overline{\partial}\omega\wedge\partial\varphi
\\
&=&\frac{1}{V_{\omega}}\int_{X}\varphi(\omega^{3}-\omega^{3}_{\varphi})
\\
&-&\frac{3}{2V_{\omega}}\int_{X}\varphi\omega\wedge\sqrt{-1}\partial\omega\wedge\overline{\partial}\varphi
+\frac{3}{2V_{\omega}}\int_{X}\varphi\omega\wedge\sqrt{-1}\overline{\partial}\omega\wedge\partial\varphi
\\
&-&\frac{3}{2V_{\omega}}\int_{X}\varphi\omega_{\varphi}\wedge\sqrt{-1}\partial\omega\wedge\overline{\partial}\varphi
+\frac{3}{2V_{\omega}}\int_{X}\varphi\omega_{\varphi}\wedge\sqrt{-1}\overline{\partial}\omega\wedge\partial\varphi.
\end{eqnarray*}
Thus

\begin{proposition} \label{C.1} One has
\begin{eqnarray*}
& & \mathcal{I}^{{\rm AY}}_{\omega}(\varphi) \ \ = \ \
\frac{1}{V_{\omega}}\int_{X}\varphi(\omega^{3}-\omega^{3}_{\varphi})\\
&-&\frac{3}{2V_{\omega}}\int_{X}\varphi\omega_{\varphi}
\wedge\sqrt{-1}\partial\omega\wedge\overline{\partial}\varphi-
\frac{3}{2V_{\omega}}\int_{X}\varphi\omega\wedge\sqrt{-1}\partial\omega\wedge
\overline{\partial}\varphi \nonumber\\
&+&\frac{3}{2V_{\omega}}\int_{X}\varphi\omega_{\varphi}
\wedge\sqrt{-1}\overline{\partial}\omega\wedge\partial\varphi
+\frac{3}{2V_{\omega}}\int_{X}\varphi\omega\wedge\sqrt{-1}\overline{\partial}\omega\wedge\partial\varphi
\nonumber.
\end{eqnarray*}
\end{proposition}

Similarly, we have
\begin{eqnarray*}
& &\mathcal{J}^{{\rm AY}}_{\omega}(\varphi) \ \ = \ \
-\mathcal{L}^{{\rm
M}}_{\omega}(\varphi)+\frac{1}{V_{\omega}}\int_{X}\varphi\omega^{3}
\\
&+&\frac{3}{8}\frac{3}{2V_{\omega}}\int_{X}\varphi\omega\wedge\sqrt{-1}\partial\omega\wedge\overline{\partial}\varphi
-\frac{3}{4}\left(2+\frac{3}{2}\right)\frac{1}{2V_{\omega}}\int_{X}\varphi\omega_{\varphi}
\wedge\sqrt{-1}\partial\omega\wedge\overline{\partial}\varphi \\
&-&\frac{3}{8}\frac{1}{2V_{\omega}}\int_{X}\varphi\omega\wedge\sqrt{-1}\overline{\partial}\omega
\wedge\partial\varphi+\frac{3}{4}\left(2+\frac{3}{2}\right)\frac{1}{2V_{\omega}}\int_{X}\varphi
\omega_{\varphi}\wedge\sqrt{-1}\overline{\partial}\omega\wedge\partial\varphi
\\
&-&\frac{3}{8}\frac{1}{2V_{\omega}}\int_{X}\varphi\omega_{\varphi}\wedge\sqrt{-1}
\partial\omega\wedge\overline{\partial}\varphi+\frac{3}{8}\frac{1}{2V_{\omega}}\int_{X}
\varphi\omega_{\varphi}\wedge\sqrt{-1}\overline{\partial}\omega\wedge\partial\varphi
\\
&-&\frac{27}{8}\frac{1}{2V_{\omega}}\int_{X}\varphi\omega\wedge\sqrt{-1}\partial\omega\wedge
\overline{\partial}\varphi+\frac{27}{8}\frac{1}{2V_{\omega}}\int_{X}\varphi\omega\wedge
\sqrt{-1}\overline{\partial}\omega\wedge\partial\varphi \\
&=&-\mathcal{L}^{{\rm
M}}_{\omega}(\varphi)+\frac{1}{V_{\omega}}\int_{X}\varphi\omega^{3}
\\
&-&\frac{3}{2V_{\omega}}\int_{X}\varphi\omega\wedge\sqrt{-1}\partial\omega\wedge
\overline{\partial}\varphi+\frac{3}{2V_{\omega}}\int_{X}\varphi\omega\wedge\sqrt{-1}
\overline{\partial}\omega\wedge\partial\varphi \\
&-&\frac{3}{2V_{\omega}}\int_{X}\varphi\omega_{\varphi}\wedge\sqrt{-1}\partial\omega\wedge\overline{\partial}\varphi
+\frac{3}{2V_{\omega}}\int_{X}\varphi\omega_{\varphi}\wedge\sqrt{-1}\overline{\partial}\omega\wedge\partial\varphi
\end{eqnarray*}

So
\begin{proposition} \label{C.2} One has
\begin{eqnarray*}
& & \mathcal{J}^{{\rm AY}}_{\omega}(\varphi) \ \ = \ \
-\mathcal{L}^{{\rm
M}}_{\omega}(\varphi)+\frac{1}{V_{\omega}}\int_{X}\varphi\omega^{3}
\\
&-&\frac{3}{2V_{\omega}}\int_{X}\varphi\omega_{\varphi}
\wedge\sqrt{-1}\partial\omega\wedge\overline{\partial}\varphi
-\frac{3}{2V_{\omega}}\int_{X}\varphi\omega\wedge\sqrt{-1}\partial\omega\wedge\overline{\partial}\varphi\nonumber\\
&+&\frac{3}{2V_{\omega}}\int_{X}\varphi\omega_{\varphi}
\wedge\sqrt{-1}\overline{\partial}\omega\wedge\partial\varphi
+\frac{3}{2V_{\omega}}\int_{X}\varphi\omega\wedge\sqrt{-1}\overline{\partial}\omega\wedge\partial\varphi.\nonumber
\end{eqnarray*}
\end{proposition}

\bibliographystyle{amsplain}

\end{document}